\documentclass[english]{article}
\usepackage[T1]{fontenc}
\usepackage[latin9]{inputenc}
\usepackage{verbatim}
\usepackage{amsmath}
\usepackage{amsthm}
\usepackage{amssymb}
\usepackage{esint}
\usepackage{varwidth}
\usepackage{url}
\usepackage{constants}
\usepackage{graphicx}
\usepackage{float}
\usepackage{hyperref}
\usepackage{mathabx} 
\usepackage{subcaption}
\usepackage{commath}
\usepackage{tikz}
\usetikzlibrary{decorations.pathreplacing}
\usepackage{authblk}
\usepackage{lipsum}
\usepackage{algorithmic}
\usepackage{algorithm}
\usepackage{placeins}
\usepackage{todonotes}

\hypersetup{
   colorlinks=true,%
   citecolor=black,%
   filecolor=black,%
   linkcolor=black,%
   urlcolor=blue
}

\makeatletter
\theoremstyle{plain}
\newtheorem{theorem}{\protect\theoremname}
  \theoremstyle{plain}
  \newtheorem{lemma}[theorem]{\protect\lemmaname}
  \theoremstyle{plain}

  \theoremstyle{plain}
  \newtheorem{definition}[theorem]{\protect\definitionname}
  \theoremstyle{plain}
  \newtheorem{proposition}[theorem]{\protect\propositionname}
  \theoremstyle{definition}
  
  \newtheorem*{remark}{Remark}
  \theoremstyle{plain}

  \theoremstyle{plain}
\newtheorem*{informal}{Informal Theorem}

\DeclareMathOperator*{\argmin}{\arg\!\min}
\DeclareMathOperator*{\argmax}{\arg\!\max}

\makeatother

\usepackage{babel}
  \providecommand{\definitionname}{Definition}
  \providecommand{\examplename}{Example}
  \providecommand{\lemmaname}{Lemma}
  \providecommand{\corrolaryname}{Corollary}
  \providecommand{\propositionname}{Proposition}
  \providecommand{\conditionsname}{Conditions}
  \providecommand{\assumptionname}{Assumption}
\providecommand{\theoremname}{Theorem}

\begin{document}

\title{The dual approach to non-negative super-resolution:\\ 
perturbation analysis and $\ell_1$ data fidelity}

\author[1,2]{St\'{e}phane Chr\'{e}tien}
\author[3]{Andrew Thompson}
\author[4]{Bogdan Toader}
\affil[1]{Laboratoire ERIC, Universit\'e Lyon 2, Bron, France} 
\affil[2]{Alan Turing Institute, London, UK}
\affil[3]{National Physical Laboratory, Teddington, UK}
\affil[4]{MRC Laboratory of Molecular Biology, Cambridge, UK}

\maketitle

\begin{abstract}
  We study the problem of super-resolution, where we 
  recover the locations and weights of non-negative point sources
  from (potentially noisy) samples of their convolution with a Gaussian kernel. 
  Previous work has shown that exact recovery is possible 
  by minimising the total variation norm of the measure, 
  and a practical way to achieve this is by solving the dual problem. 
  In this paper, we study the stability of the solution with respect 
  to the solution of the dual problem, both in the case of exact 
  measurements and in the case of measurements with additive noise.
  In particular, we establish a relationship between perturbations 
  in the dual variable and perturbations in the primal variable around the
  true solution and, in the case of inexact measurements,
  we derive a similar relationship between the additive noise and perturbations in the dual variable 
  using an $\ell_1$ data fidelity whose dual is box constrained.
  Our analysis is based on a quantitative version of the implicit function theorem.
\end{abstract}


\section{Problem setup}

In the study of non-negative super-resolution, the aim is 
to estimate a signal $x$ which consists of a number of
point sources with unknown locations and non-negative magnitudes,
from only a few measurements of the convolution of $x$ with
a known convolution kernel $\phi$. This is a problem that 
arises in a number of applications, for example 
fluorescence microscopy \cite{betzig2006imaging},
astronomy \cite{puschmann2005super}
or ultrasound imaging \cite{tur2011innovation}.
In such applications, the measurement
device has a limited resolution and cannot distinguish 
between distinct point sources that are close to each 
other in the input signal $x$. 
This is often modelled as a deconvolution problem with a Gaussian kernel.

Specifically, let $x$ be a 
non-negative measure on $I=[0,1]$ consisting
of $k$ unknown non-negative 
point sources:
\begin{equation*}
  x = \sum_{i=1}^{k} a_i \delta_{t_i},
\end{equation*}
with $a_i > 0$, for all $i=1,\ldots,k$, 
and let $y_j$ be the possibly noisy measurements
obtained by sampling 
the convolution of $x$ with a known 
kernel $\phi$ 
at locations $s_j$:
\begin{equation}
  y_j 
  = \int_I \phi(t-s_j) x(\dif t) + w_j
  = \sum_{i=1}^{k} a_i \phi(t_i-s_j) + w_j,
  \label{eq:y noise}
\end{equation} 
for all $j=1,\ldots,m$ or, in vector notation:
\begin{equation}
  y = \sum_{i=1}^k a_i \Phi(t_i) + w,
\end{equation}
where
\begin{align}
  y &= [y_1,\ldots,y_m]^T, \\
  \Phi(t) &= [\phi(t-s_1),\ldots,\phi(t-s_m)]^T, 
  \label{eq:vec phi} \\
  w &= [w_1,\ldots,w_m]^T.
\end{align}
Of particular interest is the case of the Gaussian kernel:
\begin{equation}
  \phi(t) = e^{-t^2/\sigma^2},
  \label{eq:gaussian}
\end{equation}
where $\sigma$ is assumed to be known to the practitioner.

In the setting where the measurements $y$ are exact, 
namely when $w=0$, the signal $x$ can be recovered
by solving the following problem:
\begin{equation}
  \min_{x \geq 0} \| x \|_{TV}
  \quad \text{subject to} \quad
  y = \int_I \Phi(t) x(\dif t),
  \label{eq:main prog}
\end{equation}
where $\Vert \cdot \Vert_{TV}$ is the Total Variation (TV) norm 
for Radon measures defined as
\begin{align}
    \|x\|_{\mathrm{TV}} & =\sup \left\{\int \psi \mathrm{d} x ;
      \quad \psi \in C(I),\|\psi\|_{\infty} \leqslant 1\right\}.
\end{align}
When the measurements are corrupted by additive noise, we depart from the standard approach in the super-resolution literature relying on a least squares $\ell_2$ data fidelity term, and solve instead an $\ell_1$ minimisation problem:
\begin{equation}
  \min_{x \geq 0} 
  \left\|
    y - \int_{I} \Phi(t) x(\dif t) 
  \right\|_1
  \quad \text{such that} \quad
  \| x \|_{TV} \leq \Pi,
  \label{eq:primal noisy}
\end{equation}
where $\Pi$ plays the role of the regularisation parameter.
While $\ell_1$ is also natural for robustness~\cite{pierucci2014smoothing}, our motivation is structural, as we will clarify in Section~\ref{sec:main goals} and Section~\ref{sec:pert noisy}.

In the context of problems \eqref{eq:main prog} and \eqref{eq:primal noisy},
in this manuscript we give bounds on the errors in the source
locations $\{t_i\}_{i=1}^k$ and weights $\{a_i\}_{i=1}^k$
as a function of the errors in the dual variable when solving the
dual problem, which we then extend to the case when the measurements are
corrupted by additive noise, where we give an exact dependence of the
error in the dual variable on the level of noise.
A subset of the results in this paper has been presented in the conference 
article~\cite{chretien_dual_2019}.

The problem of super-resolution has been studied extensively
in the literature since the seminal paper \cite{Candes2014}, 
which addressed the case of complex amplitudes. Since the original contributions of Cand\`es and Fernandez-Granda, there have been numerous 
follow-up results such as the ones by Schiebinger et al. \cite{schiebinger2017superresolution}, 
Duval and Peyr{\'e} \cite{duval2015exact}, Denoyelle et al. \cite{Denoyelle2017}, 
Bendory et al. \cite{bendory2016robust}, Aza\"is et al. \cite{Azais2015}, Eftekhari et al.
\cite{eftekhari_non-negative_2018, eftekhari2018superresolution}, and, more recently, Kalra et al.~\cite{kalra_small-noise_2024}, Gabet et al.~\cite{gabet_global_2025}, Poon and Peyr{\'e}~\cite{poon_super-resolved_2025}, and Carioni and Del Grande~\cite{carioni_general_2026}.
For instance, the authors of \cite{schiebinger2017superresolution} consider the noiseless setting by taking real-valued 
samples of $y$ with a more general choice of $\phi$ (such as a Gaussian) 
and also assume $x$ to be non-negative as in the present work. 
Their proposed approach again involves TV norm minimisation with linear constraints. 
Bendory et al. \cite{bendory2016robust} consider $\phi$ to be Gaussian or Cauchy, do not place sign assumptions on 
$x$, and also analyse the TV norm minimisation with linear fidelity constraints for estimating $x$ from noiseless samples of $y$.

\subsection{Motivation: dual certificate and $\ell_1$ data fidelity}\label{sec:main goals}

A standard way to approach problem \eqref{eq:main prog} is by considering its dual:
\begin{equation}
  \max_{\lambda \in \mathbb{R}^m}\;
  y^T \lambda
  \quad \text{subject to} \quad 
  \lambda^T \Phi(t) \leq 1
 \quad \forall t \in I,
  \label{eq:dual}
\end{equation}
which is a finite-dimensional problem 
with infinitely many constraints, known as a semi-infinite 
program (SIP).

Solving the dual problem for $\lambda$ leads to the \textit{dual certificate}, 
a function of the form $q(s)=\sum_{j=1}^m \lambda_j \phi(s-s_j)$ 
(defined in Section~\ref{sec:pert noise-free}),  whose
global maximisers are the source locations $\{t_i\}_{i=1}^k$. The 
weights $\{a_i\}_{i=1}^k$ are then found by solving a least squares problem
using the measurements and the source locations. Using the idea of the dual certificate,
our perturbation results are quite intuitive: the locations of the
global maximisers of the dual certificate are perturbed when $\lambda$
is perturbed, which leads to perturbed source locations $t_i$.
Providing a quantitative analysis of the recovery error 
as a function of the error in the dual solution is the first goal of the present work. 

The second goal is to extend the perturbation analysis to the case when the measurements are corrupted by additive noise. 
To do this, we leverage the specific form of the dual problem~\eqref{eq:dual} and its associated dual certificate, and consider the dual of the $\ell_1$ minimisation problem in~\eqref{eq:primal noisy}, derived in Appendix~\ref{apdx:dual noisy}:
\begin{equation}
  \max_{\substack{\beta > 0\\\lambda \in \mathbb{R}^m}}
  \beta \left(
    \lambda^T y
    - \Pi
  \right)
  \quad\text{subject to}\quad
  \lambda^T \Phi(t) \leq 1,
  \quad\forall t \in [0,1]
  \quad\text{and}\quad 
  \|\lambda\|_{\infty} \leq 1/\beta,
  \label{eq:dual noisy raw}
\end{equation}
where the box constraint is a consequence of the $\ell_1$ data term, as the convex conjugate of $\|\cdot\|_1$ is the indicator of the $\ell_{\infty}$ ball.
Since we only need the maximiser $\lambda$ of~\eqref{eq:dual noisy raw} (not the optimal value) to form the dual certificate $q(s)$, and for any fixed $\beta$ the argmax over $\lambda$ is unchanged by the positive scaling $\beta$ and the additive constant $-\beta\Pi$ in the objective, we may treat $\beta$ (equivalently the box radius $\tau = 1/\beta$) as the regularisation parameter and solve the argmax problem:
\begin{align}
  \argmax_{\lambda \in \mathbb{R}^m}
    \lambda^T y
  \quad\text{subject to}\quad
  \lambda^T \Phi(t) \leq 1,
  \quad\forall t \in [0,1]
  \quad\text{and}\quad 
  \|\lambda\|_{\infty} \leq \tau,
  \label{eq:dual noisy} 
\end{align}
where $\tau = 1/\beta$.
This specific form of the dual makes the perturbation analysis tractable in the noisy setting, as we will see in Section~\ref{sec:pert noisy}. 

Both dual problems~\eqref{eq:dual} and~\eqref{eq:dual noisy} are in the same finite-dimensional, polyhedral family, which is amenable to efficient optimisation algorithms such as exchange methods~\cite{eftekhari2019equivalence} or sequential quadratic programming~\cite{lopez2007semiinfinite}. 
Furthermore, this formulation of the dual problems can also be solved using non-smooth optimisation algorithms such as bundle methods~\cite{NesterovOpt,fan2019} by incorporating the constraint $\lambda^T \Phi(t) \leq 1, \forall t \in I$ into the objective function as a penalty term. 
For example, the noise-free dual~\eqref{eq:dual} can be reformulated as
\begin{equation}
  \min_{\lambda \in \mathbb{R}^m} 
  - y^T \lambda + C \cdot 
  \max \left\{
    \sup_{s} \left(
      \lambda^T \Phi(s) - 1
    \right),
    0
  \right\},
  \label{eq:dual penalty}
\end{equation}
for a penalty parameter $C$, with only the additional box constraint on $\lambda$ in the noisy case. 

\begin{figure}
  \centering
  \begin{subfigure}[t]{0.45\textwidth}
    \includegraphics[width=\textwidth]{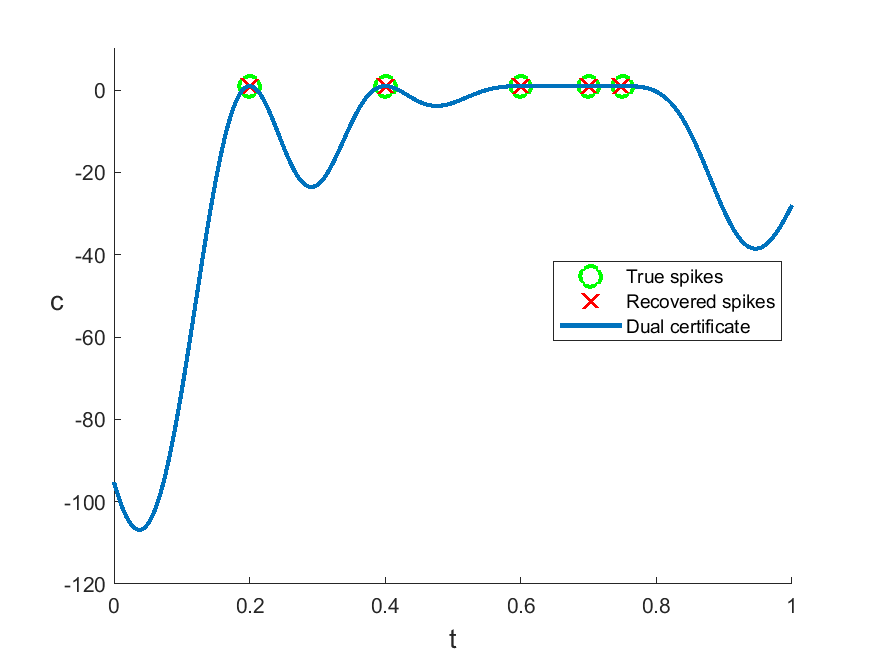}
  \end{subfigure}
  \begin{subfigure}[t]{0.45\textwidth}
    \includegraphics[width=\textwidth]{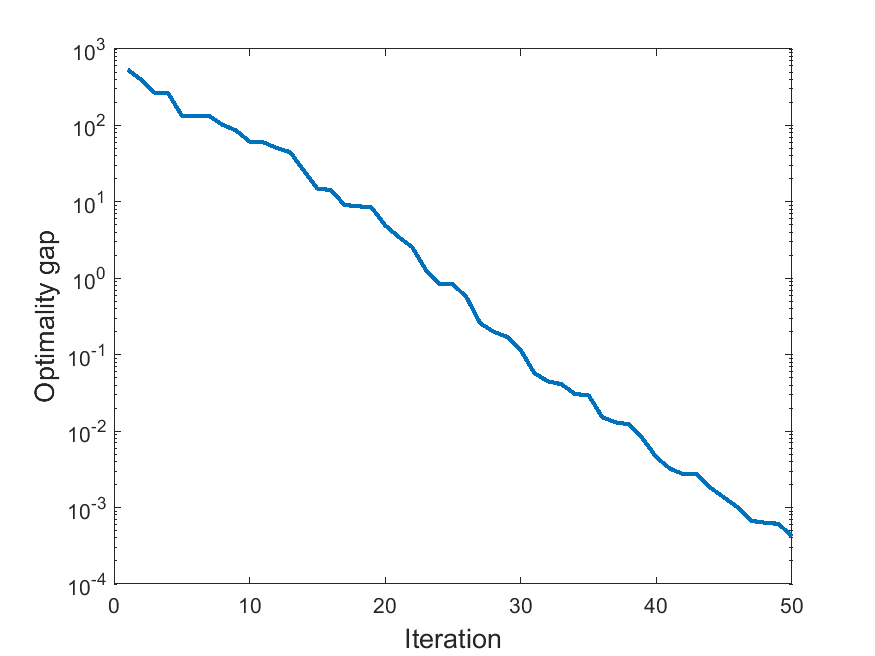}
  \end{subfigure}
  \caption{A signal consisting of five spikes with locations $\{0.2,0.4,0.6,0.7,0.75\}$ and amplitude $1$ is convolved with a Gaussian kernel $\phi(t) = e^{-t^2/\sigma^2}$ with $\sigma=0.1$ and sampled at $15$ equispaced points on $[0,1]$. Then problem~\eqref{eq:dual penalty} is solved using the level bundle method.
    (a) Computed dual certificate and solution.
    (b) Decrease in optimality gap with iteration number.}
  \label{fig:level plots}
\end{figure}

We briefly illustrate this approach by solving the noise-free dual problem~\eqref{eq:dual} via~\eqref{eq:dual penalty} using the level bundle method, introduced in~\cite{NesterovOpt} to solve continuous sparse inverse problems.
The level bundle method is applied to compute $\lambda$, and the spike locations are identified from the global maximisers of the resulting dual certificate. 
Figure~\ref{fig:level plots}(a) displays the recovered solution along with the corresponding dual certificate, showing that the method is able to recover the signal.
Figure~\ref{fig:level plots}(b) shows the speed of convergence in terms of the decrease in the optimality gap (the model gap, see Appendix~\ref{apdx:level method}), which is linear in practice.

\subsection{Contributions}

In this paper, we restrict our study to the case of Gaussian kernels. Our main results are:
\begin{itemize}
    \item When the measurements are exact, we provide bounds on how far the 
estimated spike locations $t_i$ and magnitudes $a_i$ are from
their true values as the dual variable $\lambda$ is perturbed from
the optimal solution $\lambda^*$ of~\eqref{eq:dual}. 
These bounds, given in Theorem~\ref{thm:t dep lambda} and Theorem~\ref{thm:a dep t}, offer insight into the size of the error in the spike locations and magnitudes when the super-resolution problem is solved via its convex dual.

  \item  When the measurements are corrupted by additive noise, we leverage the connection between the convex dual of~\eqref{eq:main prog} and that of the $\ell_1$ minimisation problem~\eqref{eq:primal noisy}, highlighted in Section~\ref{sec:main goals}, to study the impact of measurement noise on the solution of the noisy dual problem~\eqref{eq:dual noisy}. The main result is Theorem \ref{thm:pert noisy}, where an explicit bound is given. 
\end{itemize}

While the bounds given in these theorems apply only to the
case when the convolution kernel is Gaussian, the same
techniques can be applied to obtain perturbation bounds
for other kernels, with a few differences in the way 
some sums in the proofs are bounded, which would
be specific to the kernel used.

\subsection{Comparison with previous work}

\subsubsection{Alternative formulations for the noiseless setting}

For the particular case of non-negative $x$, Boyd et al. \cite{boyd2017} 
proposed an improved 
Frank-Wolfe algorithm in the primal. In certain instances, e.g. with Fourier 
samples (such as in \cite{Candes2014,candes2013super}), 
the dual, which is a SIP, can also be reformulated 
as a semi-definite program (SDP).
From a practical point of view, SDP is notoriously 
slow for even a moderately large number of variables.  
The algorithm of \cite{boyd2017} is a first order scheme with potential local 
correction steps, and is practically more viable.

As already mentioned, the main reason we advocate for 
using the dual problem \eqref{eq:dual} is that exact 
penalty can be used in order to reformulate the dual 
problem as a non-smooth minimisation problem for which 
methods such as bundle methods~\cite{NesterovOpt, Hiriart-Urruty1996} are efficient in practice. 
To the best of our knowledge, 
there is no 
analysis of the impact of obtaining approximate solutions 
of the dual
on the quality of the recovered locations.

\subsubsection{The penalised least squares approach}

The approach adopted in \cite{duval2015exact,Denoyelle2017} is to solve a 
least-squares-type minimisation problem with a TV norm based penalty term (also referred to as the Beurling 
LASSO; see for example \cite{de2012exact}) for recovering $x$ from samples of $y$.  
The approach in \cite{duval2017sparse} considers a natural finite approximation on the grid to the 
continuous problem, and studies the limiting behaviour as the grid becomes finer; see also \cite{duval2017sparse2}.
These works develop a perturbation analysis which is different from ours since it applies to specific types of perturbations of a different problem ($\ell_2$ vs. $\ell_1$ type fidelity terms), and does not provide precise quantitative dependencies with respect to all the parameters of the problem.
Choosing the $\ell_1$ fidelity enables the analysis in the noisy case presented in Section~\ref{sec:pert noisy}.

\subsubsection{The Prony/Matrix Pencil approach}

Another efficient approach is the one of~\cite{Moitra2015} based on the original work of Hua and Sarkar~\cite{hua1990} using a Matrix Pencil approach, extended to the multi-kernel setting in~\cite{Chretien2020}, and further investigated in~\cite{katz_decimated_2023, liu_super-resolution_2023}. Perturbation analysis of the Matrix Pencil approach is provided in \cite{Moitra2015}; see also \cite{Chretien2020} for a more detailed exposition of these results with the correct order of dependencies. The reason we develop an analysis of the dual problem \eqref{eq:dual} here is that it easily extends to the multidimensional setting as well, at least for small dimensions. In contrast, the Matrix Pencil method, although very efficient in one dimension, becomes much more involved in several dimensions~\cite{josz2019sparse}. 

\subsection{Plan of the paper}

We start by presenting
 the noise-free perturbation results 
related to problem \eqref{eq:dual} in Section \ref{sec:pert noise-free},
followed by the perturbation results in the
setting where the measurements are corrupted by noise
in Section \ref{sec:pert noisy}. The proofs of 
our results are given in Section \ref{sec:proofs} 
and we show numerical experiments to verify
the validity of our results in practice
in Section \ref{sec:lvl numerics}.
Lastly, we
conclude the paper in Section \ref{sec:conclusion}.

\section{Bound on the error as $\lambda$ is perturbed -- the noise-free case}
\label{sec:pert noise-free}

In this section we present our first main results, namely
two theorems that give bounds on the perturbations
around the source locations $t_i$ and the 
magnitudes $a_i$ respectively, as the dual variable
is perturbed away from the optimiser $\lambda^*$, 
when the convolution kernel is a 
Gaussian with known width $\sigma$ 
as defined in \eqref{eq:gaussian}.

First, let us briefly give an informal statement of the
main results in this section.

\begin{informal}
  \label{thm:primal informal}
  \textbf{(Stability of primal recovery)}
  Let $\lambda^* \in \mathbb{R}^m$ be a solution of the dual 
  program \eqref{eq:dual} with $\phi$ Gaussian and $\lambda$
  a perturbation of $\lambda^*$ in a ball of radius $\delta_{\lambda}$,
  given in Theorem~\ref{thm:t dep lambda}
  and let $\mathbf{t}^*, \mathbf{a^*}$ be the vectors of source
  locations and weights in the true signal $x$ 
  and $\mathbf{\tilde{t}},\mathbf{\tilde{a}}$ respectively
  their perturbations due to $\lambda$. 
  Then, the error between $\mathbf{t^*}$ and $\tilde{\mathbf{t}}$
  is bounded by:
  \begin{equation}
    \|\mathbf{\tilde{t}}-\mathbf{t^*}\|_2
      \leq \sqrt{k} C_{t^*} \|\lambda-\lambda^*\|_2.
  \end{equation}
  Moreover, if the above error is bounded 
  by $\delta_t$, given in Theorem~\ref{thm:a dep t}, then
  the error between $\mathbf{a^*}$ and $\tilde{\mathbf{a^*}}$
  is bounded by:
  \begin{equation}
    \|\mathbf{\tilde{a}}-\mathbf{a^*}\|_2 
      \leq
      C_{a^*} 
      \|\tilde{\mathbf{t}}-\mathbf{t^*}\|_2
      + O(\|\tilde{\mathbf{t}}-\mathbf{t^*}\|_2^2).
  \end{equation}
\end{informal}

As the two error bounds above are derived independently
using different ideas, we will discuss them individually.
Before giving the exact statement of each theorem, 
we define the concept of dual certificate,
which plays an important role throughout this paper.
\begin{definition}
  \textbf{(Dual certificate)}
  Consider a solution $\lambda^*$ of the dual 
  problem \eqref{eq:dual}  or \eqref{eq:dual noisy}.
  Then a dual certificate is a function of the form
  \begin{equation}
    q(t) = \sum_{j=1}^{m} \lambda_j^* \phi(t-s_j)
    = {\lambda^*}^T \Phi(t),
    \label{eq:dual cert}
  \end{equation}
  which satisfies the conditions:
  \begin{align}
    &q(t_i) = 1, \quad \forall i=1,\ldots,k,
    \label{eq:cond1}
    \\
    &q(t) < 1, \quad \forall t \ne t_i, 
        \forall i=1,\ldots,k.
    \label{eq:cond2}
  \end{align}
  \label{def:dual certif}
\end{definition}
The idea of dual certificate is common in the super-resolution
literature, and we know that 
the global maximisers of $q(t)$ correspond to the 
source locations $\{t_i\}_{i=1}^k$ (see, for 
example~\cite{Candes2014,schiebinger2017superresolution,eftekhari2018superresolution},
Once these are found,
amplitudes $\{a_i\}_{i=1}^k$ are obtained
by solving a linear system.

We are now ready to discuss the perturbation results 
in the noise-free setting. In the following theorem,
we consider 
the dual \eqref{eq:dual} of \eqref{eq:main prog} and 
quantify how the source locations given by the global
maximisers of the dual certificate formed by the
dual solution $\lambda^*$ are affected by 
perturbations of $\lambda^*$.
\begin{theorem}
  \label{thm:t dep lambda}
  \textbf{(Dependence of $|t-t^*|$ on $\|\lambda-\lambda^*\|_2$)}
  Let $\lambda^* \in \mathbb{R}^m$ be a solution of the 
  dual program \eqref{eq:dual} with $\phi$ Gaussian 
  as given in \eqref{eq:gaussian} such that the the dual
  certificate $q(s)$ defined in \eqref{eq:dual cert}
  satisfies conditions \eqref{eq:cond1} and \eqref{eq:cond2},
  $\lambda$ a perturbation of $\lambda^*$ in a ball of 
  radius $\delta_{\lambda}$ and $t$ an arbitrary local
  maximiser of $q_{\lambda}(s)=\sum_{j=1}^m \lambda_j \phi(s-s_j)$.
  Note that, for $\lambda=\lambda^*$, 
  the corresponding local (and global) maximiser $t^*$ 
  of $q_{\lambda^*} = q$
  is a true source location in $\{t_i\}_{i=1}^k$.
  Let $R = \frac{ \| \lambda^* \|_2}{\sigma}$ and
  $c \approx 3.9036$ a universal constant. If
  the radius $\delta_{\lambda}$ is bounded by
  \begin{equation}
    \delta_{\lambda} \leq
    \frac{
      |q''(t^*)|^2 \sigma^3 \sqrt{e}
    }{
      4\sqrt{2} \left(
        2 + c 
        R 
      \right)
      m
    },
    \label{eq:delta_lambda}
  \end{equation}
  then
  \begin{equation}
    |t-t^*| \leq C_{t^*} \|\lambda-\lambda^*\|_2,
    \label{eq:t dep lambda}
  \end{equation}
  where
  \begin{align}
    C_{t^*} &= \frac{1}{4 + cR}
      \left[
        1 + \frac{
          2\sqrt{2m} (2+cR)
        }{
          |q''(t^*)| \sqrt{e} 
        }
      \right].
    \label{eq:C_t}
    \\
  \end{align}
\end{theorem}

\begin{proposition} 
  \label{prop:simplified ct}
  \textbf{(Simplified $C_{t^*}$)}
  Under the conditions of Theorem~\ref{thm:t dep lambda},
  the constant $C_{t^*}$ can be further bounded by:
  \begin{equation}
    C_{t^*} < \frac14 + \frac{2\sqrt{2}}{\sqrt{e}}
      \cdot \frac{\sqrt{m}}{|q''(t^*)|}.
  \end{equation}
\end{proposition}

The proofs of Theorem~\ref{thm:t dep lambda} 
and Proposition~\ref{prop:simplified ct} are given 
in Section~\ref{sec:proof thm t dep lambda}.
As a brief summary, Theorem~\ref{thm:t dep lambda}
is proved by applying the implicit
function theorem to the function $F(t,\lambda) = q'(t)$, where
$q(t)$ is the dual certificate given in Definition~\ref{def:dual certif},
since we know that $F(t^*,\lambda^*) = 0$. This allows us to express
$t$ as a function $t(\lambda)$ in a neighbourhood of $(t^*,\lambda^*)$,
and a quantitative version of the theorem~\cite{bahsoun_mean_2025}
gives an explicit expression for $\partial_{\lambda}t(\lambda)$ and the 
neighbourhood in terms of the derivatives of $F$. 
By bounding this derivative and the neighbourhood  and then
applying a truncated Taylor expansion to $t(\lambda)$, we
obtain the result of Theorem~\ref{thm:t dep lambda}.

One of the main conclusions which can be drawn from this result 
is that the primal spike location error is controlled 
in $l_\infty$, but degrades as a function of the number of measurements 
in the order of $\sqrt{m}$. 
Alternatively, we can write \eqref{eq:t dep lambda} in 
terms of the $\ell_2$ norm of the error between the vector of 
true source locations $\mathbf{t^*}$
and the perturbed source locations
$\mathbf{\tilde{t}}$:
\begin{equation*}
  \|\mathbf{\tilde{t}}-\mathbf{t^*}\|_2
  \leq \sqrt{k} C_{t^*} \|\lambda-\lambda^*\|_2.
\end{equation*}

Of crucial importance is the curvature 
of the dual certificate at the true solution: the flatter the 
certificate, the worse the estimation error. 
Our theorem also gives important information about the accuracy 
in the dual variable required to guarantee our upper bound on the 
error of recovery. This accuracy is of the inverse order of the 
number of measurements, which is quite a stringent constraint. 
Both the $m$ and the $\sqrt{m}$ factors are a consequence of the
way we bound sums of shifted copies of the kernel, namely
$\sum_{j=1}^m \phi(t-s_j) \leq m \max_{t \in \mathbb{R}} \phi(t)$.
Given the fast decay of the Gaussian, it is clear that this is not
a tight bound. However, any bound would reflect the density
of samples close to each source location.

We will now give a result regarding the perturbation
of the magnitudes $a_i$ when $\lambda^*$ is perturbed.
Let $\Phi$ be the matrix whose entries are defined as
\begin{equation}
    \Phi_{ij} = \phi(t_j-s_i),
    \label{eq:def Phi}
\end{equation}
and $\mathbf{t^*}$ and $\mathbf{a^*}$ the vectors of 
source locations and weights:
\begin{align*}   
    \mathbf{t^*} = [t_1,\ldots,t_k]^T,
    \quad\quad
    \mathbf{a^*} = [a_1,\ldots,a_k]^T.
\end{align*}
When we solve \eqref{eq:dual} exactly, we obtain
the source locations by finding the global maximisers
of $q(s)$. Then, the vector of weights $\mathbf{a^*}$ 
is found by solving the system
\begin{equation*}
    \Phi \mathbf{a} = \mathbf{y}.
\end{equation*}
When the source locations are perturbed, 
we denote the resulting perturbed data matrix by:
\begin{equation}
  \tilde{\Phi} = \Phi + E,
  \label{eq:phi tilde def}
\end{equation}
and we calculate the vector of perturbed weights
$\tilde{\mathbf{a}}$ as the solution of the least 
squares problem
\begin{equation}
  \min_{\mathbf{a}} \| \tilde{\Phi} \mathbf{a} - \mathbf{y} \|_2.
  \label{eq:perturbed phi least sq}
\end{equation}
The following theorem, proved in 
Section \ref{sec:proof thm a dep t},
gives a bound on the error 
$\| \mathbf{a^*} - \tilde{\mathbf{a}} \|_2$ 
between the vector of true weights $\mathbf{a^*}$
and the vector of weights $\tilde{\mathbf{a}}$ obtained by solving
the least squares problem  \eqref{eq:perturbed phi least sq}
with the perturbed matrix $\tilde{\Phi}$, as a function
of the error $\|\tilde{\mathbf{t}}-\mathbf{t^*}\|_2$ 
between the perturbed source 
locations $\tilde{\mathbf{t}}$ and the true source 
locations $\mathbf{t^*}$.

\begin{theorem}
  \label{thm:a dep t}
  \textbf{(Dependence of $\|\mathbf{\tilde{a}} - \mathbf{a^*}\|_2$ 
    on $\|\mathbf{\tilde{t}}-\mathbf{t^*}\|_2$)}
  Let $\mathbf{t^*} \in [0,1]^{k}$ be the vector of true source 
  locations, $\tilde{\mathbf{t}} \in [0,1]^k$ the perturbed source 
  locations in a ball of radius $\delta_t$, 
  $\mathbf{a^*}$ the vector of true weights 
  and $\tilde{\mathbf{a}}$ the vector of perturbed weights obtained by
  solving problem \eqref{eq:perturbed phi least sq}.
  If the radius $\delta_t$ is bounded by:
  \begin{equation}
    \delta_t
    <
    \frac{
        \sigma^2 \sigma_{\max}(\Phi)
    }{
        4e^{4/\sigma^2} \sqrt{m}
    }
    \left(
        \sqrt{
            1 + \frac{
                \sigma_{\min}^2(\Phi)
            }{
                \sigma_{\max}^2(\Phi)
            } 
        } -1
    \right),
    \label{eq:cond t}
  \end{equation} 
  where $\sigma_{\max}(\Phi)$, $\sigma_{\min}(\Phi)$
  are the largest and respectively smallest singular values
  of the matrix $\Phi$ defined in \eqref{eq:def Phi},
  then:
  \begin{equation}
    \|\mathbf{\tilde{a}}-\mathbf{a^*}\|_2 
    \leq
    C_{a^*} 
    e^{\frac{4}{\sigma^2}\max_j |\tilde{t}_j-t_j|}
    \|\tilde{\mathbf{t}}-\mathbf{t^*}\|_2
    + O(\|\tilde{\mathbf{t}}-\mathbf{t^*}\|_2^2),
    \label{eq:a dep t}
  \end{equation}
  where
  \begin{equation*}
    C_{a^*} = 
    \frac{
      4 \sqrt{m}
      \|\mathbf{a^*}\|_2
    }{
      \sigma^2 \sigma_{\min}(\Phi)
    }.
  \end{equation*}
\end{theorem}

Note that we write the $O(\|\tilde{\mathbf{t}}-\mathbf{t^*}\|_2^2)$
term in the bound above in order to simplify the 
presentation of the result.
We can, however, calculate the constants 
corresponding to the higher order terms in the bound
by using the inequality \eqref{eq:a bound expansion} 
in the proof of Theorem \ref{thm:a dep t} in
Section \ref{sec:proof thm a dep t}. For example,
the constant in the second order term is equal to
$  
  C_{a^*}^2/\|\mathbf{a^*}\|_2
  \left[
    1 + 2 \sigma^2_{\max}(\Phi) /
    \sigma^2_{\min}(\Phi)
  \right]
$.

\section{Bound on $\|\lambda-\lambda^*\|_2$
    in terms of the noise $w$}
\label{sec:pert noisy}

In this section we assume that the measurements
are corrupted by additive noise
and we give a result where
we bound the perturbation in the dual variable $\lambda$
around the minimiser $\lambda^*$ as a function of the
noise $w$ in the measurements. Specifically,
the noisy measurements are defined as in \eqref{eq:y noise}:
\begin{equation*}
  y_j = \int_I \phi_j(t) x(\dif t) + w_j
  = \sum_{i=1}^{k} a_i \phi_j(t_i) + w_j,
\end{equation*}
for $w_j \ne 0$ and $j=1,\ldots,m$.

The aim is to estimate how 
the source locations $\{t_i\}_{i=1}^k$
and weights $\{a_i\}_{i=1}^k$ are affected by 
the additive noise $w$ in the measurements
around the solution of the problem.
In the previous section, we have established how
the source locations and weights are perturbed
around their true values 
as the dual variable $\lambda$ is perturbed
around its optimal value $\lambda^*$. 
In the noisy setting, we establish a quantitative relationship
between the perturbations of $\lambda$ around $\lambda^*$
and the magnitude of the noise.

Before we state the main result, we first describe the exact mathematical 
setting under which the result holds. Then, we introduce 
the function $\bar{F}$ in \eqref{eq:F bar} to which we apply the 
implicit function theorem, whose Jacobian is crucial for this result.

As discussed in Section~\ref{sec:main goals}, we consider the $\ell_1$ minimisation problem~\eqref{eq:primal noisy} to account for noise in the measurements, and we focus on a particular form of its dual problem, given in~\eqref{eq:dual noisy}, which we restate here:
\begin{align*}
  \argmax_{\lambda \in \mathbb{R}^m}
  y^T \lambda 
  \quad &\text{such that} \quad
  \lambda^T \Phi(t) \leq 1,
    \quad \forall t \in I,
  \nonumber \\
  &\quad\quad \text{and} \quad
  \| \lambda \|_{\infty} \leq \tau,
\end{align*}
The parameter $\tau$ is the inverse of the Lagrange multiplier
corresponding to the constraint in \eqref{eq:primal noisy},
and therefore it plays the same regularisation role as $\Pi$
(see Section~\ref{sec:main goals} and Appendix~\ref{apdx:dual noisy}).

To motivate the exact form of the function $\bar{F}$ 
in \eqref{eq:F bar} to which we apply the implicit function
theorem to obtain the perturbation result from 
Theorem \ref{thm:pert noisy}, consider the exact penalty
formulation of \eqref{eq:dual noisy}:
\begin{equation}
  \argmin_{\lambda \in \mathbb{R}^m} \Psi_{\Pi}(\lambda)
  \quad \text{such that} \quad
  \| \lambda \|_{\infty} \leq \tau,
  \label{eq:exact pen}
\end{equation}
where
\begin{equation}
  \Psi_{\Pi}(\lambda) =
  - y^T \lambda + \Pi \cdot 
  \max \left\{
    \sup_{s} \left(
      \sum_{j=1}^m \lambda_j 
        \phi(s-s_j) - 1
    \right),
    0
  \right\}.
\end{equation}
For a large enough value of $\Pi$, a solution
to \eqref{eq:exact pen} which satisfies the 
constraints in \eqref{eq:dual noisy} is also
a solution of \eqref{eq:dual noisy} 
(see, for example, Section 1.2 in \cite{Hiriart-Urruty1996}).
This is a non-smooth optimisation problem and its solution
can be found by using any method that relies on calculating
subgradients, for example the level method \cite{NesterovOpt}.

A subgradient of $\Psi_{\Pi}(\lambda)$ has the form:
\begin{equation}
  \partial \Psi_{\Pi} =
  \begin{cases}
    -y + \Pi \sum_{i=1}^{k'} 
      \nu_i g(s_i^*),
      \quad (\nu_1+\ldots+\nu_{k'} = 1)
      &\quad \text{ if } \quad
      \sup_{s} 
        \sum_{j=1}^m \lambda_j \phi(s-s_j) 
      > 1,
    \\
    -y + \Pi \sum_{i=1}^{k'} 
      \nu_i g(s_i^*),
      \quad (\nu_1+\ldots+\nu_{k'} \leq 1)
      &\quad \text{ if } \quad
      \sup_{s} 
        \sum_{j=1}^m \lambda_j \phi(s-s_j) 
      = 1,
    \\
    -y,
      &\quad \text{ if } \quad
      \sup_{s} 
        \sum_{j=1}^m \lambda_j \phi(s-s_j) 
      < 1,
  \end{cases}
\end{equation}
where $\{s_i^*\}_{i=1}^{k'}$ are the global
maximisers of the function
$
  \sum_{j=1}^m \lambda_j \phi(s-s_j) 
$,
the vectors $g(s)$ are of the form
$
  g(s) = [\phi(s-s_1),\ldots,\phi(s-s_m)]^T
$
and $\nu_i \geq 0$ for all $i=1,\ldots,k'$.
Note that here we apply the formula 
for the subgradient of the $\max$ function and for 
the $\sup$ function (see for example \cite{Hiriart-Urruty1996}). 
The coefficients in the convex 
combination from the formula for the subgradient
of the $\max$ function with zero account for the 
case when $\nu_1 + \ldots + \nu_{k'} < 1$\footnote{
  More specifically, both functions in the $\max$ 
  attain their maximum, so we have that
  \newline
  $
    \partial \Psi_{\pi} = 
    -y + \Pi 
    \left[ \alpha_1 \partial 
      \sup_s \left(
        \sum_{j=1}^m \lambda_j \phi(s-s_j)-1
      \right)
    + \alpha_2 \partial 0 \right]
  $,
  with $\alpha_1,\alpha_2>0$ 
  and $\alpha_1+\alpha_2=1$,
  and therefore
  $
    \partial \Psi_{\pi} = 
    -y + \Pi 
    \sum_{i=1}^{k'} \alpha_1 \nu_i'g(s_i^*)
  $,
  with $\nu_1'+\ldots+\nu_{k'}'=1$ 
  and $0 \leq \alpha_1 \leq 1$.
}.

As in the noise-free setting, we assume here that
the dual solution $\lambda^*$ forms a dual certificate,
namely the function $q(s)$ as defined in \eqref{eq:dual cert}
satisfies conditions \eqref{eq:cond1} and \eqref{eq:cond2}.
Then, the subdifferential at $\lambda^*$ has the form:
\begin{equation}
  \partial \Psi_{\Pi}(\lambda^*) = 
  -y + \Pi \sum_{i=1}^k \nu_i g(t_i),
\end{equation}
where $\{t_i\}_{i=1}^k$ are the source
locations, so the optimality condition 
for \eqref{eq:exact pen}:
\begin{equation}
  0 \in \partial \Psi_{\Pi}(\lambda^*),
\end{equation}
is equivalent to:
\begin{equation}
  y = \Pi \sum_{i=1}^k \nu_i g(t_i),
  \label{eq:opt cond}
\end{equation}
for some $\nu_1,\ldots,\nu_k \geq 0$
with $\nu_1 + \ldots + \nu_k \leq 1$
and for $w=0$.
Note that, given the definition 
of $y$ from \eqref{eq:y noise},
the optimality condition \eqref{eq:opt cond}
is satisfied for
\begin{align}
  \nu_i &= \frac{a_i}{\Pi},
    \quad \forall i=1,\ldots,k, 
  \\
  w_j &= 0, 
    \quad \forall j=1,\ldots,m,
\end{align}
where in order to satisfy $\nu_1+\ldots+\nu_k \leq 1$,
we need $\Pi$ such that:
\begin{equation}
  \Pi \geq a_1+\ldots+a_i,
\end{equation}
which is the same as the constraint in \eqref{eq:primal noisy}.

Motivated by the above reasoning, we now apply 
the quantitative implicit function theorem,
as given in~\cite{bahsoun_mean_2025},
to a function $F$ of the form:
\begin{equation}
  F([\lambda,\nu]^T, w) = 
  \sum_{i=1}^k a_i \Phi(t_i^*) - 
  \sum_{i=1}^k \nu_i \Phi(t_i(\lambda)) + w,
  \label{eq:F QIFT noisy}
\end{equation}
where we know that $F([\lambda^*,a]^T, 0) = 0$.
For the sake of simplicity, we include 
the parameter $\Pi$ in the coefficients $\nu_i$,
so in the second sum in $F$ each $\nu_i$ 
actually corresponds to $\Pi \nu_i$, 
and $\nu_1+\ldots+\nu_k \leq \Pi$ 
rather than $\nu_1+\ldots+\nu_k \leq 1$.

However, note that 
$
  F: \mathbb{R}^{m+k} \times \mathbb{R}^{m}
  \to 
  \mathbb{R}^m
$
and in order to apply the implicit function theorem
to obtain the dependence of the first argument of $F$
as a function of the second argument, it is required
that the spaces of the first argument, the second argument
and the codomain of $F$ have the same dimension.
To overcome this issue, we assume that the solution 
we work with has a few particular properties,
since the dual certificate,
given in Definition~\ref{def:dual certif}, is not unique
in general. As before, we will assume that the 
solution $\lambda^*$ to the dual problem \eqref{eq:dual noisy}
satisfies the dual certificate condition. In addition, 
we assume the existence of a solution $\bar{\lambda}^*$
of \eqref{eq:dual noisy} as follows:
\begin{definition} 
  Let $\lambda^* \in \mathbb{R}^m$ be a solution to the dual 
  problem \eqref{eq:dual noisy} with $m-k$ entries on the
  boundary of the box constraint of \eqref{eq:dual noisy}
  i.e. there exist indices 
  $\gamma_1,\ldots,\gamma_{m-k} \in \{1,\ldots,m\}$ such that
  $\lambda^*_{\gamma_j} = \pm \tau$, $j = 1,\ldots,m-k$.
  Then we define $\bar{\lambda}^* \in \mathbb{R}^{k}$ to
  be the vector that consists of the non-fixed
  entries of $\lambda^*$, in the same order,
  and $\bar{\lambda} \in \mathbb{R}^{k}$ a perturbation
  of $\bar{\lambda}^*$.
  \label{def:lambda bar}
\end{definition}

\begin{remark}
  In practice, such a solution $\bar{\lambda}^*$ would be achieved due to the complementarity conditions at optimality corresponding to the box constraint $\|\lambda\| \leq \tau$. 
  This is a direct consequence of using an $\ell_1$ data fidelity term in~\eqref{eq:primal noisy}: an $\ell_2$ fidelity term would instead give a Euclidean ball constraint (or a smooth quadratic dual term) and none of this complementarity.
\end{remark}

Similarly, we define a vector consisting of $2k$ entries
of $\Phi(t)$ in~\eqref{eq:vec phi}.
\begin{definition}
  Let $\theta := \{\theta_1,\ldots,\theta_{2k}\} 
  \subset \{1,\ldots,m\} $. 
  Then we define $\bar{\Phi}_{\theta}(t)$ to be
  the vector consisting of the entries of $\Phi(t)$
  in \eqref{eq:vec phi} corresponding to the indices
  in $\theta$:
  \begin{equation}
    \bar{\Phi}_{\theta}(t) := [
      \phi(t-s_{\theta_1}),\phi(t-s_{\theta_2}),
      \ldots, \phi(t-s_{\theta_{2k}})
    ]^T.
  \end{equation}
  We will also use $\bar{\Phi}(t)$ to denote 
  $\bar{\Phi}_{\theta}(t)$ when the specific choice
  of $\theta$ is not relevant in the context.
  \label{def:vec phi bar}
\end{definition}  
Lastly, given the definitions of $\bar{\lambda}$ 
and $\bar{\Phi}_{\theta}(t)$ above,
we define the following function to which we will
be able to apply the implicit function theorem:
\begin{definition}
  Let $\lambda^* \in \mathbb{R}^m$ be a 
  solution of \eqref{eq:dual noisy}
  with $m-k$ fixed entries and $\bar{\lambda}^* \in \mathbb{R}^k$ 
  consisting of 
  the non-fixed entries of $\lambda^*$, as given 
  in Definition~\ref{def:lambda bar}, 
  and let $\bar{\Phi}_{\theta}(t)$ be given as in 
  Definition~\ref{def:vec phi bar} for an index
  set $\theta$ of $2k$ indices between $1$ and $m$.
  Then, for 
  the perturbation $\bar{\lambda}$ of $\bar{\lambda}^*$,
  we define the function
  $
    \bar{F}: \mathbb{R}^{2k} \times \mathbb{R}^{2k}
    \to 
    \mathbb{R}^{2k}
  $ 
  as: 
  \begin{equation}
    \bar{F}([\bar{\lambda},\nu]^T, w_{\theta}) = 
    \sum_{i=1}^k a_i \bar{\Phi}_{\theta}(t_i^*) - 
    \sum_{i=1}^k \nu_i \bar{\Phi}_{\theta}(t_i(\bar{\lambda})) + w_{\theta},
    \label{eq:F bar}
  \end{equation}
  where $t_i(\bar{\lambda}^*) = t_i^*$, for $i=1,\ldots,k$,
  are the source locations corresponding to $\lambda^*$ 
  and $w_{\theta} \in \mathbb{R}^{2k}$ contains the entries of
  the noise vector $w \in \mathbb{R}^m$ corresponding to the
  indices in $\theta$.
  \label{def:F bar}
\end{definition}

We can now state the main result of this section, namely
a bound on the perturbation of $\lambda^*$ (or more specifically
$\bar{\lambda}^*$) as a function of the measurement noise.
The proof is given 
in Section \ref{sec:pert noisy proof}. 

\begin{theorem}
  \label{thm:pert noisy}
  \textbf{(Dependence of $\|\lambda-\lambda^*\|_2$ on the noise $w$)}
  Let $\lambda^* \in \mathbb{R}^m$ be a solution to the 
  dual problem \eqref{eq:dual noisy} with $w=0$ ,
  namely the optimal solution of \eqref{eq:dual noisy}
  with noiseless measurements, which satisfies the 
  conditions in Definition~\ref{def:lambda bar}, and the 
  vector $\bar{\lambda}^* \in \mathbb{R}^{k}$ of non-fixed
  entries of $\lambda^*$.
  For the function $\bar{F}$ in Definition~\ref{def:F bar},
  let $J^*$ be its Jacobian with respect 
  to the first variable, evaluated 
  at $([\bar{\lambda}^*, a]^T, 0)$ and $\sigma_{\min}(J^*)$ its smallest
  singular value. We also assume that the solution $\lambda^*$ 
  forms a dual
  certificate, namely the function $q(t)$ defined 
  in \eqref{eq:dual cert} satisfies 
  conditions \eqref{eq:cond1} and \eqref{eq:cond2}.
  If $J^*$ is invertible, $\|w\|_2 \leq \delta_w$ and
  \begin{equation}
    |q''(t^*)| \leq 2 \left(
    1 + \frac{4 m \tau}{\sigma^2} 
  \right),
  \label{eq:cond q thm}
  \end{equation}
  then, for a perturbation $\lambda$ of $\lambda^*$ with 
  the same fixed entries to the boundary of the box 
  constraint, we have that:
  \begin{equation}
    \left\|
      \lambda - \lambda^*
    \right\|_2 \leq
    C_{\lambda^*} \cdot
    \left\|
      w
    \right\|_2,
    \label{eq:pert noisy}
  \end{equation}
  where
  \begin{align}
    &C_{\lambda^*} = \frac{2}{\sigma_{\min}(J^*)},
    \\
    &\delta_{w} =
    \frac{\sigma_{\min}(J^*)^2}{
      4P(m,k,\sigma,\Pi,\tau,C_{t^*})
    },
  \end{align}
  and
  \begin{align}
    P(&m,k,\sigma,\Pi,\tau,C_{t^*}) 
      =\sqrt{2} k
      \Bigg[
        \frac{1}{\sigma^2} 
        \left(
          2 \sqrt{k} C_{t^*}^2 \Pi
          + 4 k C_{t^*} \bar{\Delta}_2 \tau \Pi
        \right)
      \nonumber \\
      &+ \frac{1}{\sigma} 
      \left(
        \frac{\sqrt{2k} C_{t^*}}{\sqrt{e}}
        + 4 \sqrt{k} C_{t^*}^2 \Pi
        + \frac{2\sqrt{2} \bar{\Delta}_2 \Pi}{\sqrt{e}}
        + 8 k C_{t^*} \bar{\Delta}_2 \tau \Pi
        + \frac{\sqrt{2k} \bar{\Delta}_2 \Pi}{\sqrt{e}}
        + \sqrt{\frac2e} C_{t^*}
      \right)
      \Bigg],
    \label{eq:P def}   
  \end{align}
  where $C_{t^*}$ is given in \eqref{eq:C_t}
  in Theorem~\ref{thm:t dep lambda},
  $c \approx 3.9036$, 
  $c_2 = 4 + \frac{c\sqrt{2}}{\sqrt{e}} \approx 7.3484$
  are universal constants and 
  \begin{equation}
    \bar{\Delta}_2 =
    \frac{\sqrt{k}}{\sigma^4} \left(
      c_2 C_{t^*} m \tau
    + \frac{2\sqrt{2}}{\sqrt{e}} \sigma
    \right).
    \label{eq:delta 2 bar} 
  \end{equation}
\end{theorem}

The theorem above makes explicit the dependence of the
perturbation in the dual variable $\lambda$ around
the solution $\lambda^*$ on the additive noise $w$ 
in the measurement vector $y$, with the assumption given
in Definition~\ref{def:lambda bar}.
This is a linear relation
where the constant depends on the specific configuration
of the problem we are solving, namely the locations 
and weights of the sources, and width of the Gaussian and
the sampling locations.
The theorem also gives an upper bound on the magnitude
of the noise where this result holds as a function
of the same parameters.

As an additional interpretation of Theorem~\ref{thm:pert noisy}
regarding the assumption on the fixed entries in $\lambda$
and $\lambda^*$, it
states that, for a solution $\lambda$ to the dual
problem \eqref{eq:dual noisy} with noisy measurements
that has $m-k$ entries equal to the boundary of the box 
constraint, there is a solution $\lambda^*$ to the 
noise-free dual problem $w=0$ with the same entries fixed
to the boundary of the box constraint and the error for
the remaining $k$ entries bounded 
by \eqref{eq:pert noisy}.

Moreover, under a few additional assumptions, we give a 
simplified approximation of the constant $P$ in \eqref{eq:P def}
for clarity:
\begin{proposition}
  \label{prop:simplified p}
  \textbf{(Simplified P)}
  Under the conditions of Theorem~\ref{thm:pert noisy} and,
  in addition, if $\Pi,\tau \leq 1$, then:
  \begin{equation}
    P(m,k,\sigma,C_{t^*}) 
    = \mathcal{O}\left(
      \frac{m k^{5/2} C_{t^*}^2}{\sigma^6}
    \right).
  \end{equation}
\end{proposition}

One important observation is that, while the above result
only applies to a subset of the entries in $\lambda$ and $w$,
which entries are selected is not arbitrary. 
The choice of the entries in $\lambda$ and $w$ reflects
which samples $s_j$ contain the most information, and therefore
which noise entries in $w$ affect the solution to the 
optimisation problem the most. More specifically,
in order for the Jacobian $J^*$ to be invertible,
we are led to select the samples (and therefore $\lambda$ and $w$
entries) that satisfy this condition the best, namely the
ones that are the closest to the source locations. We discuss 
this aspect in more detail in Section \ref{sec:discussion}.

Lastly, note that the results in Section \ref{sec:pert noise-free}
and Section \ref{sec:pert noisy} refer to different optimisation
problems: the duals \eqref{eq:dual} and \eqref{eq:dual noisy}
of problems \eqref{eq:main prog} and \eqref{eq:primal noisy} 
respectively. However, the proofs of our perturbation results rely on
the property that the dual solution $\lambda$ forms a dual
certificate, the global maximisers of which give the 
locations of the point sources in the input signal $x$,
with the additional bound on $\lambda$ from \eqref{eq:dual noisy}
being used in the proof of Theorem \ref{thm:pert noisy}.
Moreover, since our analysis is independent of the exact
formulation of the primal problems, we conclude that
the results from both Section \ref{sec:pert noise-free} 
and Section \ref{sec:pert noisy} apply to the 
problem of super-resolution in the noisy setting, 
namely they give bounds of the perturbations 
of the source locations and weights as a consequence
of noise in the measurements.

\subsection{Discussion}
\label{sec:discussion}

One of the conditions in Theorem~\ref{thm:pert noisy}
is that the Jacobian $J^*$ is invertible. While we do not
provide a rigorous analysis of the conditions in which 
this is satisfied, in this section we discuss in more 
detail what the condition requires and give further motivation
for why it is true in a reasonable scenario.
Specifically, we assume that the samples that are used for
calculating the Jacobian are the closest samples to the sources,
i.e. the set $\theta$ for which we define $\bar{F}$
in Definition~\ref{def:F bar} contains 
the two indices corresponding to the closest two samples 
to each source location, for each of the $k$ sources.
Therefore, the rows in the system given by $\bar{F}$ 
in \eqref{eq:F bar}, as well as 
the entries in $\bar{\lambda}$ and the entries
in the noise vector $w_{\theta}$,
correspond to these samples.

Recall that $J^*$ is the Jacobian of the
function $\bar{F}$ from \eqref{eq:F bar}
with respect to the first argument. 
The entries in $J^*$ are:
\begin{align}
  \partial_{\bar{\lambda}_l} \bar{F}([\bar{\lambda},\nu]^T,w) 
    \Bigr|_{\substack{\bar{\lambda}=\bar{\lambda}^*\\\nu=a\\ w_{\theta}=0}}
  &= -\sum_{i=1}^k a_i \phi'(t_i^* - s_{\theta_j}) 
    \partial_{\bar{\lambda}_l} t_i(\bar{\lambda}^*) \\
  &=
    \sum_{i=1}^k \frac{
      a_i \phi'(t_i^*-s_{\theta_j}) \phi'(t_i^*-s_l)
    }{
      q''(t_i^*)
    },
    \label{eq:jac block left}
\end{align}
for $l=1,\ldots,k$, $j=1,\ldots,2k$,
where $\{s_l\}_{l=1}^{k}$ correspond to the 
non-fixed entries of $\lambda$ (i.e. $\bar{\lambda}$)
and
\begin{equation}
  \partial_{\nu_l} \bar{F}([\bar{\lambda},\nu]^T,w_{\theta})
  \Bigr|_{\substack{\bar{\lambda}=\bar{\lambda}^*\\\nu=a\\ w_{\theta}=0}}
  = - \phi(t_l^*-s_{\theta_j}), 
  \label{eq:jac block right}
\end{equation}
for $l=1,\ldots,k$, $j=1,\ldots,2k$,
where in the first equality we used \eqref{eq:qift res}
with \eqref{eq:partial t} and \eqref{eq:partial lambda} 
plugged in, so the result holds under the conditions 
in Theorem~\ref{thm:t dep lambda}, namely 
for $\lambda$ with $\|\lambda-\lambda^*\|_2 \leq \delta_{\lambda}$,
where $\delta_{\lambda}$ is given 
in \eqref{eq:delta_lambda}.

Writing $J^*$ as 
\begin{equation}
  J^* = [J_{\lambda} J_{\nu}],
\end{equation}
where the entries in the blocks $J_{\lambda}$
and $J_{\nu}$ are given by \eqref{eq:jac block left}
and \eqref{eq:jac block right} respectively, 
we have that:
\begin{equation}
  J_{\lambda} = \sum_{i=1}^k
    \frac{a_i}{q''(t_i^*)} \bar{\Phi}'(t_i^*) \bar{\Phi}'(t_i^*)^T,
\end{equation}
and 
\begin{equation}
  J_{\nu} = - [\bar{\Phi}(t_1^*) \ldots \bar{\Phi}(t_k^*)],
\end{equation}
where
\begin{align}
  &\bar{\Phi}(t) = [\phi(t-s_{\theta_1}),\ldots,\phi(t-s_{\theta_{2k}})]^T,
  \\
  &\bar{\Phi}'(t) = [\phi'(t-s_{\theta_1}),\ldots,\phi'(t-s_{\theta_{2k}})]^T.
\end{align}
Note that $rank(J_{\nu}) = k$ by the T-systems property of 
the Gaussian (assuming that the $t_1 \leq \ldots \leq t_k$
and $s_{\theta_1} \leq \ldots \leq s_{\theta_{2k}}$) and in order for the 
matrix $J^*$ to be invertible we need $rank(J^*)=2k$,
as it is a square matrix with $2k$ columns. 
By rewriting the columns of $J_{\lambda}$, we have that:
\begin{align}
  J^* = \left[
    \sum_{i=1}^k \frac{a_i \phi'(t_i^* - s_1)}{q''(t_i^*)} \bar{\Phi}'(t_i^*) 
    \quad \ldots \quad
    \sum_{i=1}^k \frac{a_i \phi'(t_i^* - s_k)}{q''(t_i^*)} \bar{\Phi}'(t_i^*) 
    \quad
    - \bar{\Phi}(t_1^*)
    \quad \ldots \quad
    - \bar{\Phi}(t_k^*)
  \right],
\end{align}
and by taking its determinant and using the multi-linearity property
of the determinant with respect to its columns, we have that:
\begin{align}
  \det(J^*) &= (-1)^k
  \frac{
    a_1\ldots a_k
  }{
    q''(t_1^*) \ldots q''(t_k^*)
  } 
  \nonumber \\
  &\cdot \sum_{l=1}^{k!} 
    \left(
      \prod_{i=1}^k \phi'(P_l(t_i^*)-s_i)
    \right)
    \left|
      P_l(\bar{\Phi}'(t_1^*)) 
      \quad \ldots \quad
      P_l(\bar{\Phi}'(t_k^*))
      \quad
      \bar{\Phi}(t_1^*)
      \quad \ldots \quad
      \bar{\Phi}(t_k^*)
    \right|,
  \label{eq:det J expansion}
\end{align}
where $P_l$ for $l=1,\ldots,k!$ are the permutations 
of $k$ elements. 
Note that when we expand the determinant, the terms in the final sum
are determinants with all the possible combinations of the vectors in each
sum, which results in most determinants having repeated columns, so they are
equal to zero. The only non-zero determinants in the resulting sum are the ones
where the first $k$ columns are the vectors $\{\Phi'(t_i^*)\}_{i=1}^k$ 
and their permutations, multiplied by the corresponding constants.
We now order the columns of the determinant:
\begin{align}
  \det(J^*) &= (-1)^k
  \frac{
    a_1\ldots a_k
  }{
    q''(t_1^*) \ldots q''(t_k^*)
  } 
  \nonumber \\
  &\qquad\qquad \sum_{l=1}^{k!} 
    sign(P_i)
    \left(
      \prod_{i=1}^k \phi'(P_l(t_i^*)-s_i)
    \right)
    \left|
      \bar{\Phi}(t_1^*)
      \quad
      \bar{\Phi}'(t_1^*)
      \quad \ldots \quad
      \bar{\Phi}(t_k^*)
      \quad
      \bar{\Phi}'(t_k^*)
    \right|,
    \nonumber \\
  &= (-1)^k 
  \frac{
    a_1\ldots a_k
  }{
    q''(t_1^*) \ldots q''(t_k^*)
  } 
  \left|
    \bar{\Phi}(t_1^*)
    \quad
    \bar{\Phi}'(t_1^*)
    \quad \ldots \quad
    \bar{\Phi}(t_k^*)
    \quad
    \bar{\Phi}'(t_k^*)
  \right|
  \nonumber \\
  &\qquad\qquad \sum_{l=1}^{k!} 
    sign(P_l)
    \left(
      \prod_{i=1}^k \phi'(P_l(t_i^*)-s_i)
    \right),
\end{align}
where by $sign(P_i)$ we denote the sign of the determinant
corresponding to the permutation $P_i$ after reordering 
the columns as above.
Because of the extended T-system property of the Gaussian
function~\cite{karlin1966tchebycheff}, the determinant
above is strictly positive. 
The dominant term in the sum is the one corresponding to the
identity permutation,  where for each $i=1,\ldots,k$, the sample $s_i$ is
the closest sample to the source location $t_i^*$.
As the samples get further, the terms of the sum approach zero.
This can be expressed more quantitatively by imposing
explicit conditions on the distances between the closest samples and the sources,
the separation of sources and the separation of samples, 
as done, for example, in \cite{eftekhari2018superresolution}.

As a last remark to motivate the choice of the dimension 
of $\bar{\lambda}^*$ and $\bar{\lambda}$ in Definition~\ref{def:lambda bar}, 
note the expansion in \eqref{eq:det J expansion}
of $\det(J^*)$.
If the vector $\bar{\lambda}$ had more than $k$ entries,
then the columns consisting 
of the permutations of $\bar{\Phi}'(t_i^*)$ would inevitably be
repeated, since there are $k$ sources $t_i^*$ and more than $k$
such columns. This implies that all the determinants in the
sum would be zero and, therefore, $J^*$ would not be invertible,
implying that Theorem~\ref{thm:pert noisy} would not be true 
in this case. This explains why choosing $\bar{\lambda}$ to 
contain more than $k$ entries of $\lambda$ would be incompatible
with our analysis in the proof of Theorem~\ref{thm:pert noisy}.

\section{Proofs}
\label{sec:proofs}

In this section we present the proofs of the theorems 
from Sections \ref{sec:pert noise-free} and \ref{sec:pert noisy}.

\subsection{Proof of Theorem~\ref{thm:t dep lambda}
  (Dependence of $|t-t^*|$ on $\|\lambda-\lambda^*\|_2$)}
  \label{sec:proof thm t dep lambda}

Let $t^*$ be an arbitrary local maximiser of the 
function $q(t)$ in \eqref{eq:dual cert}, so $t^*$ 
is also a source location, and $\lambda^*$ 
the solution to \eqref{eq:dual}.
The key step in this proof is applying a 
quantitative version of the 
implicit function theorem~\cite{bahsoun_mean_2025}
to the function:
\begin{equation}\label{eq:Fdef}
  F(t,\lambda) = \sum_{j=1}^{m} \lambda_j \phi'(t-s_j),
\end{equation}
where $F(t^*,\lambda^*) = 0$ because $t^*$ is a maximizer of $q(s)$ in \eqref{eq:dual cert}. 
The theorem allows us to express $t$ 
as a function $t(\lambda)$ of $\lambda$ with:
\begin{equation}
    \partial_{\lambda} t(\lambda) = 
    -\left[
        \partial_t F(t(\lambda),\lambda)
    \right]^{-1}
    \partial_{\lambda} F(t(\lambda),\lambda),
    \label{eq:qift res}
\end{equation}
for $t$ in a ball of radius $\delta_0$ around $t^*$ and 
for $\lambda$ in a ball of radius $\delta_1 \leq \delta_0$ 
around $\lambda^*$, where $\delta_0$ is chosen such that
\begin{equation}
  \sup_{(t,\lambda) \in V_{\delta}}
  \left\| 
    I - 
    \left[
      \partial_t F(t^*,\lambda^*)
    \right]^{-1}
    \partial_t F(t,\lambda)
  \right\|
  \leq \frac12,
  \label{eq:cond delta0}
\end{equation}
where
$
  V_{\delta}= \left\{ 
    (t, \lambda) \in \mathbb{R}^{m+1}:
    |t - t^*| \leq \delta_0,
    \|\lambda - \lambda^*\| \leq \delta_0
  \right\}
$  
and $\delta_1$ is given by 
\begin{equation}
    \delta_1 = (2M_t B_{\lambda})^{-1} \delta_0,
    \label{eq:cond delta1}
\end{equation}
where
\begin{align*}
  B_{\lambda} &=
  \sup_{(t,\lambda) \in V_{\delta}}
    \| \partial_{\lambda} F(t,\lambda) \|_2,
  \\
  M_t &= \left\| 
    \partial_t F(t^*, \lambda^*)^{-1} 
  \right\|_2.
\end{align*}
The following two lemmas, proved in
Sections~\ref{sec:proof Lemma bnd t} and \ref{sec:proof Lemma bnd lambda}
respectively, give us values of $\delta_0$
and $\delta_1$ that define balls around $t^*$ and $\lambda^*$
respectively which are included in the balls required
by the quantitative implicit function theorem with radii 
defined in \eqref{eq:cond delta0} and \eqref{eq:cond delta1}.
\begin{lemma}
  \textbf{(Radius of ball around $t^*$)}
  The condition \eqref{eq:cond delta0} is satisfied if
  \begin{equation}
      \delta_0 
      =\frac{
          \sigma^2 |q''(t^*)| 
        }{
          \sqrt{m}\left(
            4 + 2c \cdot
            \frac{\|\lambda^*\|_2}{\sigma}
          \right)
        }.
      \label{eq:apriori delta bound}
  \end{equation}
  \label{lem:bnd t}
\end{lemma}

\begin{lemma}
  \textbf{(Radius of ball around $\lambda^*$)}
  For $\delta_0$ from Lemma~\ref{lem:bnd t} 
  and $\delta_1$ from condition \eqref{eq:cond delta1}, 
  the following choice of $\delta_{\lambda}$: 
  \begin{equation*}
    \delta_{\lambda} 
    = \frac{
      \sigma \sqrt{e} |q''(t^*)| 
    }{
      2\sqrt{2m}
    } \cdot \delta_0
  \end{equation*}
  satisfies $\delta_{\lambda} < \delta_1$.
  \label{lem:bnd lambda}
\end{lemma}
Given the definition of the function $F$ in \eqref{eq:Fdef}, we have that
\begin{align}
  \partial_{t} F(t,\lambda) &= \sum_{j=1}^{m} \lambda_j \phi''(t-s_j), 
  \label{eq:partial t}
  \\
  \partial_{\lambda} F(t,\lambda) &= [\phi'(t-s_1), \quad \ldots \quad, \phi'(t-s_m)]^T. 
  \label{eq:partial lambda}  
\end{align}
By applying Taylor expansion to $t(\lambda)$ around
$\lambda^*$ in the region defined by $\delta_0$ and $\delta_{\lambda}$,
we have that
\begin{equation*}
  t(\lambda) = t(\lambda^*) + 
  \left<
    \lambda-\lambda^*,
    \partial_{\lambda} t(\lambda_{\delta})
  \right>,
\end{equation*}
for some $\lambda_{\delta}$ on the line segment 
determined by $\lambda^*$ and $\lambda$,
so
\begin{align}
  \left|
    t(\lambda)-t(\lambda^*)
  \right|  
  &\leq
  \left\|
    \lambda-\lambda^*
  \right\|_2
  \cdot
  \left\|
    \partial_{\lambda} t(\lambda_{\delta})
  \right\|_2
  \nonumber \\
  &\leq
  \frac{\delta_0}{
    \sum_{j=1}^m {\lambda_{\delta}}_j \phi''(t(\lambda_{\delta})-s_j)
  }
  \cdot \left\| \left[
    \phi'(t(\lambda_{\delta})-s_1),
    \quad \ldots \quad, 
    \phi'(t(\lambda_{\delta})-s_m)
  \right] \right\|_2,
  \label{eq:terms_bnd}
\end{align}
where in the last inequality we used that
$\|\lambda - \lambda^*\| \leq \delta_0$ and \eqref{eq:qift res}.
We now need to bound the terms in \eqref{eq:terms_bnd} for the 
Gaussian kernel $\phi(t) = e^{-t^2/\sigma^2}$.
First, we rewrite the last inequality as
\begin{align}
  |t(\lambda) - t(\lambda^*)| 
  &\sum_{j=1}^m ({\lambda_{\delta}}_j + \lambda_j^* - \lambda_j^*)
    \phi''(t(\lambda_{\delta})-s_j)
  \nonumber \\
  &\leq
  \delta_0
  \cdot \left\| \left[
    \phi'(t(\lambda_{\delta})-s_1),
    \quad \ldots \quad, 
    \phi'(t(\lambda_{\delta})-s_m)
  \right] \right\|_2,
\end{align}
we apply the reverse triangle inequality in the sum on the left hand side:
\begin{align}
  |t(\lambda) - t(\lambda^*)| 
  &\left[
    -\left|
      \sum_{j=1}^m ({\lambda_{\delta}}_j - \lambda_j^*)
        \phi''(t(\lambda_{\delta})-s_j)
    \right| 
    + \left|
      \sum_{j=1}^m \lambda_j^* \phi''(t(\lambda_{\delta})-s_j)
    \right|
  \right]
  \nonumber \\
  &\leq
  \delta_0
  \cdot \left\| \left[
    \phi'(t(\lambda_{\delta})-s_j)
  \right]_{j=1}^m \right\|_2
\end{align}
and then we apply the Cauchy-Schwartz inequality to the first sum 
on the left hand side above to obtain:
\begin{align}
  |t(\lambda) - t(\lambda^*)| 
  &\left[
    -\left\| 
      \lambda_{\delta} - \lambda^*
    \right\|_2
    \cdot
    \left\|
      \left[ 
        \phi''(t(\lambda_{\delta})-s_j)
      \right]_{j=1}^m
    \right\|_2 
    + \left|
      \sum_{j=1}^m \lambda_j^* \phi''(t(\lambda_{\delta})-s_j)
    \right|
  \right]
  \nonumber \\
  &\leq
  \delta_0
  \cdot \left\| \left[
    \phi'(t(\lambda_{\delta})-s_j)
  \right]_{j=1}^m \right\|_2.
\end{align}
To simplify the notation, we write
$\delta_t = |t(\lambda)-t(\lambda^*)|$ and
\begin{align}
  A &=
    \left\|
      \left[ 
        \phi''(t(\lambda_{\delta})-s_j)
      \right]_{j=1}^m
    \right\|_2,
  \\
  B &=
    \left|
      \sum_{j=1}^m 
      \lambda_j^* \phi''(t(\lambda_{\delta})-s_j)
    \right|,
    \label{eq:B def}
  \\
  C &= 
    \left\| \left[
      \phi'(t(\lambda_{\delta})-s_j)
    \right]_{j=1}^m \right\|_2,
\end{align}
and by using\footnote{
    Since $\|\lambda - \lambda^*\|\leq \delta_0$
    and $\lambda_{\delta}$ is on the line segment
    between $\lambda^*$ and $\lambda$, then
    $\lambda_{\delta}$ is in the ball centred at $\lambda^*$
    with radius $\delta_0$.
} 
$\|\lambda_{\delta}-\lambda^*\|_2 \leq \delta_0$,
we have that:
\begin{equation}
  \delta_t
  (-\delta_0 A + B)
  \leq
  \delta_0 C,
\end{equation}
which can be further re-written as:
\begin{equation}
  \delta_t \leq
  \frac{C + \delta_t A}{B}
  \cdot \delta_0.
  \label{eq:equiv terms bnd}
\end{equation}
The aim now is to obtain a bound on $\delta_t$ as a function
of $\delta_0$ and the parameters of the problem. 
Therefore, we need to lower bound $B$ and 
upper bound $C+\delta_t A$.

\subsubsection*{Bounding A,B,C}

We start with $B$, for which we want to calculate a lower bound.
First, we Taylor expand each term of the sum around
$ t(\lambda^*) - s_j $ as follows:
\begin{align}
  B &= \left|
      \sum_{j=1}^m
      \lambda_j^* \phi''(t(\lambda^*)-s_j + t(\lambda_{\delta}) - t(\lambda^*))
  \right|
  \nonumber \\
  &= \left|
      \sum_{j=1}^m
      \lambda_j^* \phi''(t(\lambda^*)-s_j)
      + \left(t(\lambda_{\delta}) - t(\lambda^*)\right)
      \sum_{j=1}^m
      \lambda_j^* \phi'''(\xi_j)
  \right|
  \\
  &\geq \left|
      \sum_{j=1}^m
      \lambda_j^* \phi''(t(\lambda^*)-s_j)
    \right| - \left|
      t(\lambda_{\delta}) - t(\lambda^*)
    \right| \left|
      \sum_{j=1}^m
      \lambda_j^* \phi'''(\xi_j)
  \right|,
\end{align}
where 
$
  \xi_j \in 
  \left[
    t(\lambda^*) - s_j - |t(\lambda_{\delta})-t(\lambda^*)|,\quad
    t(\lambda^*) - s_j + |t(\lambda_{\delta})-t(\lambda^*)|
  \right]
$ for $j = 1, \ldots, m$, and on the last line we 
used the reverse triangle inequality. 
We calculate an upper bound of the last sum in the previous equation as follows:
\begin{align}
  \left|
    \sum_{j=1}^m
    \lambda_j^* \phi'''(\xi_j)
  \right|
  &\leq
  \left\| 
    \lambda^* 
  \right\|_2 \cdot
  \left\| \left[
    \phi'''(\xi_j)
  \right]_{j=1}^m \right\|_2,
  \quad \quad \text{by Cauchy-Schwartz,}
  \\
  &\leq
  \frac{
    c \| \lambda^* \|_2 \sqrt{m}
  }{
    \sigma^3
  },
  \label{eq:bound lambda phi'''}
\end{align}
where in the last line we used the maximum value
of $\phi'''(t)$ and $c$ is a constant.\footnote{
  \label{fn:derivs}
  $\max_{t \in \mathbb{R}} \phi'(t) 
    = \frac{\sqrt{2}}{\sigma\sqrt{e}},
  \max_{t \in \mathbb{R}} \phi''(t) 
    = \frac{2}{\sigma^2},
  \max_{t \in \mathbb{R}} \phi'''(t) 
    = \frac{c}{\sigma^3},
  \text{ where } 
  c=\frac{4\sqrt{9-3\sqrt{6}}}{e^{\frac{3-\sqrt{6}}{2}}} \approx 3.9036.$
}

Finally, by using 
the $\delta_0$ from Lemma~\ref{lem:bnd t}
as a bound on $|t(\lambda_{\delta}) - t(\lambda^*)|$ and
\eqref{eq:bound lambda phi'''}, we obtain:
\begin{equation}
  B \geq |q''(t^*)|
  \left[
    1 - \frac{
      c \| \lambda^* \|_2 
    }{
      4 \sigma + 2c \| \lambda^* \|_2
    }
  \right].
    \label{eq:bound B}
\end{equation}
Note that the last fraction above is subunitary, so the bound is indeed positive.

Lastly, we upper bound $C + \delta_t A$. We bound both $A$ and $C$ using 
the upper bounds on $\phi'$ and $\phi''$ 
given in footnote \ref{fn:derivs} and obtain:
\begin{align}
  A &\leq \frac{2\sqrt{m}}{\sigma^2},
  \label{eq:bound A}
  \\
  C &\leq \frac{\sqrt{2m}}{\sigma\sqrt{e}},
  \label{eq:bound C}
\end{align}
and for $\delta_t$ we use the bound \eqref{eq:apriori delta bound}. 
Putting \eqref{eq:apriori delta bound}, \eqref{eq:bound B}, \eqref{eq:bound A}
and \eqref{eq:bound C} together, 
we obtain:
\begin{equation}
  \left|
    t(\lambda) - t(\lambda^*)
  \right|
  \leq
  C_{t^*} \cdot 
  \left\|
    \lambda - \lambda^*
  \right\|_2,
\end{equation}
where
\begin{equation}
  C_{t^*} = \frac{
    2 \sqrt{2m} 
    \left(
      2 \sigma + c \| \lambda^* \|_2
    \right)
  }{
    | q''(t^*) | \sigma \sqrt{e}
    \left(
      4\sigma + c \| \lambda^* \|_2
    \right)
  } + \frac{
    2 \sigma 
  }{
    4 \sigma + F \| \lambda^* \|_2 
  },
  \label{eq:ct def}
\end{equation}
which can also be written in the form 
in \eqref{eq:C_t} in Theorem~\ref{thm:t dep lambda}.

\subsubsection{Proof of Lemma~\ref{lem:bnd t}  
    (Radius $\delta_0$ of the ball around 
    $t^*$)}
  \label{sec:proof Lemma bnd t}

Let us now find the radius $\delta_0$ 
which satisfies \eqref{eq:cond delta0}. Using \eqref{eq:partial t}, 
the expression inside the $\sup$ in \eqref{eq:cond delta0} is
\begin{equation}
  E = 
  \left|
    1 - \frac{
      \sum_{j=1}^m \lambda_j \phi''(t-s_j)
    }{
      \sum_{j=1}^m \lambda_j^* \phi''(t^*-s_j)
    }
  \right| 
  = 
  \frac{
    \left|
      \sum_{j=1}^m 
        \lambda_j^* \phi''(t^*-s_j) - \lambda_j \phi''(t-s_j)
    \right|
  }{
    \left|
      \sum_{j=1}^m \lambda_j^* \phi''(t^*-s_j)
    \right|.
  }
\end{equation}
By denoting each term in the sum in the numerator in the last equation above
by $T_j$ and then adding and subtracting $\lambda_j^*$ and $t^*$, we obtain:
\begin{align}
  T_j &= \lambda_j^* \phi''(t^*-s_j)
  - (\lambda_j - \lambda_j^*)\phi''(t-s_j) 
  - \lambda_j^* \phi''(t^*-s_j+ t - t^*) 
  \nonumber \\
  &=
  -(\lambda_j - \lambda_j^*) \phi''(t-s_j) 
  -\lambda_j^* (t-t^*) \phi'''(\xi_j),
\end{align}
for some 
$
  \xi_j \in \left[
    t^*-s_j-|t-t^*|, t^*-s_j+|t-t^*|
  \right]
$. Then:
\begin{align}
  E &\leq
  \frac{
    \left|
      \sum_{j=1}^m (\lambda_j-\lambda_j^*)
      \phi''(t-s_j)
    \right| + \left|
      \sum_{j=1}^m \lambda_j^*(t-t^*)
      \phi'''(\xi_j)
    \right|
  }{
    \left|
      \sum_{j=1}^m \lambda_j^* \phi''(t^*-s_j)
    \right|
  }
  \nonumber \\
  &\leq
  \frac{
    \left\|
      \lambda-\lambda^*
    \right\|_2
    \left\|\left[
      \phi''(t-s_j)
    \right]_{j=1}^m 
    \right\|_2
    + 
    \left| t - t^* \right|
    \left|
      \sum_{j=1}^m \lambda_j^* 
      \phi'''(\xi_j)
    \right|
  }{
    \left|
      \sum_{j=1}^m \lambda_j^* \phi''(t^*-s_j)
    \right|
  } = E'
\end{align}
We now have that
\begin{align}
  \sup_{(t,\lambda) \in V_{\delta_0}} E
  &\leq
  \sup_{\substack{|t - t^*| \leq \delta_0,\\
    \|\lambda - \lambda^*\| \leq \delta_0}
  } E'
  \\
  &\leq
  \delta_0 \cdot 
  \frac{
    \left\|\left[
      \phi''(t-s_j)
    \right]_{j=1}^m 
    \right\|_2
    + 
    \left|
      \sum_{j=1}^m \lambda_j^* 
      \phi'''(\xi_j)
    \right|
  }{
    \left|
      \sum_{j=1}^m \lambda_j^* \phi''(t^*-s_j)
    \right|
  }.
\end{align}
We now further upper bound the fraction on the last line of the previous equation.
The terms in the numerator are bounded by taking the maxima of 
the functions $\phi''$ and $\phi'''$ from 
footnote \ref{fn:derivs} respectively:
\begin{equation}
  \left\|\left[
      \phi''(t-s_j)
    \right]_{j=1}^m 
  \right\|_2
  = \sqrt{\sum_{j=1}^m \phi''(t-s_j)^2}
  \leq \sqrt{m \cdot \max_j |\phi''(t-s_j)|^2}
  \leq \frac{2\sqrt{m}}{\sigma^2}
\end{equation}
and
\begin{align}
  \left|
    \sum_{j=1}^m \lambda_j^* 
    \phi'''(\xi_j)
  \right|
  &\leq \| \lambda^* \|_2
  \left\|
    \left[
      \phi'''(\xi_j)
    \right]_{j=1}^m
  \right\|_2
  \quad\quad \text{by Cauchy-Schwartz}
  \\
  &= \| \lambda^* \|_2 
  \sqrt{\sum_{j=1}^m \phi'''(\xi_j)^2}
  \\
  &\leq \| \lambda^*\|_2 
  \max_j |\phi'''(\xi_j)| \sqrt{m} 
  \\
  &= c \cdot \frac{
    \|\lambda^*\|_2 \sqrt{m}
  }{\sigma^3}, 
  \label{eq:above def c}
\end{align}
where $c = \frac{4\sqrt{9-3\sqrt{6}}}{e^{\frac{3-\sqrt{6}}{2}}} \approx 3.9036$. 
By writing 
\begin{equation}
  q(t) = \sum_{j=1}^m \lambda_j^* \phi(t-s_j)
  \label{eq:q notation}
\end{equation}
and using 
the above bounds, we have that
\begin{equation}
  \sup_{(t,\lambda) \in V_{\delta_0}} E
  \leq \delta_0 \cdot 
  \frac{
    \frac{2\sqrt{m}}{\sigma^2}
    +
    c \cdot \frac{
      \|\lambda^*\|_2 \sqrt{m}
    }{\sigma^3}
  }{
    \left| q''(t^*) \right| 
  }
  \label{eq:almost delta bound}
\end{equation}
Finally, in order to satisfy condition \eqref{eq:cond delta0}, 
we need to impose the condition that the right hand side 
of \eqref{eq:almost delta bound} is less than or equal to $\frac12$. 
We select $\delta_0$ to be the largest value that satisfies this, so:
\begin{equation}
  |t - t^*| \leq
  \delta_0
  =\frac{
    \left| q''(t^*) \right| 
  }{
    \frac{4\sqrt{m}}{\sigma^2}
    +
    2c \cdot \frac{
      \|\lambda^*\|_2 \sqrt{m}
    }{\sigma^3}
  }
  =\frac{
    \sigma^2 |q''(t^*)| 
  }{
    \sqrt{m}\left(
      4 + 2c \cdot
      \frac{\|\lambda^*\|_2}{\sigma}
    \right)
  }.
\end{equation}

\subsubsection{Proof of Lemma~\ref{lem:bnd lambda}
  (Radius $\delta_{\lambda}$ of the ball 
  around $\lambda^*$)}
  \label{sec:proof Lemma bnd lambda}

The radius $\delta_{\lambda}$ of the perturbation 
of $\lambda^*$ is given by:
\begin{equation}
  \delta_{\lambda} = (2M_t B_{\lambda})^{-1} \delta_0,
\end{equation}
where
\begin{align}
  B_{\lambda} &=
  \sup_{(t,\lambda) \in V_{\delta}}
    \| \partial_{\lambda} F(t,\lambda) \|_2,
  \\
  M_t &= \left\| 
    \partial_t F(t^*, \lambda^*)^{-1} 
  \right\|_2.
\end{align}
For $B_{\lambda}$, we have:
\begin{align}
  \| \partial_{\lambda} F(t,\lambda) \|_2
  =
  \sqrt{\sum_{j=1}^m \phi'(t-s_j)^2}
  \leq
  \frac{\sqrt{2m}}{\sigma\sqrt{e}},
\end{align}
where we have used the global maximum of the first derivative
of the Gaussian from footnote \ref{fn:derivs}, so by 
taking $\sup$ on both sides in the last equation, 
we obtain:
\begin{equation}
  B_{\lambda} \leq \frac{\sqrt{2m}}{\sigma\sqrt{e}}.
\end{equation}
Note that here we do not use any assumptions on the locations
of the sources $t_i$ and the samples $s_j$. If we did, we would
be able to obtain a tighter bound than by only using the 
absolute maximum of the function.

For $M_t$, note that we have
\begin{equation}
  M_t = |q''(t^*)|^{-1}, 
\end{equation}
where $q(t)$ is defined in \eqref{eq:q notation},
so 
\begin{equation}
  (2M_t B_{\lambda})^{-1}\delta_0 \geq
  \frac{
    \sigma \sqrt{e} |q''(t^*)| 
  }{
    2\sqrt{2m}
  } \cdot \delta_0,
\end{equation}
We then take $\delta_{\lambda}$ to be equal to the lower
bound in the equation above:
\begin{equation}
  \delta_{\lambda} = 
  \frac{
    \sigma \sqrt{e} |q''(t^*)| 
  }{
    2\sqrt{2m}
  } \cdot \delta_0,
\end{equation}
and, after substituting our choice of $\delta_0$
from \eqref{eq:apriori delta bound}, 
we obtain the radius \eqref{eq:delta_lambda}
in Theorem~\ref{thm:t dep lambda}.

\subsubsection{Proof of Proposition~\ref{prop:simplified ct}}
  \label{sec:proof prop simplified ct}

Starting from the definition of $C_{t^*}$ in \eqref{eq:C_t},
we have that:
\begin{align}
  C_{t^*} &= \frac{1}{4 + cR}
    \left[
      1 + \frac{
        2\sqrt{2m} (2+cR)
      }{
        |q''(t^*)| \sqrt{e} 
      }
    \right] 
  \nonumber \\
  &< 
  \frac{1}{4 + c \|\lambda^*\|_2/\sigma}
  + \frac{2\sqrt{2m}}{|q''(t^*)|\sqrt{e}}
  \nonumber \\
  &<
  \frac{1}{4} + \frac{2\sqrt{2}}{\sqrt{e}} 
  \cdot \frac{\sqrt{m}}{|q''(t^*)|},
\end{align}
where in the first inequality we used the definition
of $R = \frac{\|\lambda^*\|_2}{\sigma}$ 
and $\frac{2 + cR}{4 + cR} < 1$ and in the second
inequality we used $c\|\lambda^*\|_2/\sigma > 0$,
where $c \approx 3.9036$ is a universal constant.

\subsection{Proof of Theorem~\ref{thm:a dep t}
  (Dependence of $\|\mathbf{\tilde{a}} - \mathbf{a^*}\|_2$ 
    on $\|\mathbf{\tilde{t}}-\mathbf{t^*}\|_2$)}
  \label{sec:proof thm a dep t}

We apply equation (4.2) in \cite{Stewart90perturbationtheory},
with $e = 0$ (the noise in the observations), 
and obtain
\begin{equation}
  \mathbf{\tilde{a}} = \mathbf{a^*} 
    - \Phi^{\dagger} E\mathbf{a^*} - F^T E \mathbf{a^*},
  \label{eq:a tilde}
\end{equation}
where $\Phi^{\dagger} = (\Phi^T \Phi)^{-1} \Phi^T$ is the
pseudo-inverse of $\Phi$ and $F=O(E)$ is the perturbation
of the $\Phi^{\dagger}$ due to the perturbation $E$ 
of $\Phi$, namely
\begin{equation*}
  \tilde{\Phi}^{\dagger} = \Phi^{\dagger} + F^T.
\end{equation*}
In order to obtain an explicit expression for $F$, 
we write $\tilde{\Phi}^{\dagger}$:
\begin{align}
  \tilde{\Phi}^{\dagger} 
  &= (\tilde{\Phi}^T \tilde{\Phi})^{-1} \tilde{\Phi}^T 
  \nonumber \\
  &= \left[ (\Phi+E)^T (\Phi+E) \right]^{-1} (\Phi+E)^T
  \quad\quad\text{by } \eqref{eq:phi tilde def}
  \nonumber \\
  &= (\Phi^T \Phi + \Delta)^{-1} (\Phi^T + E^T),
  \label{eq:tilde phi expl}
\end{align}
where 
\begin{equation}
  \Delta = E^T \Phi + \Phi^T E + E^T E
  \in \mathbb{R}^{k \times k}.
  \label{eq:delta def}
\end{equation}
In order to compute the first factor in \eqref{eq:tilde phi expl},
consider the QR decomposition of $\Phi$:
\begin{equation}
  \Phi = QR, 
  \quad \text{where} \quad
  Q \in \mathbb{R}^{m \times k} 
  \quad \text{and} \quad
  R \in \mathbb{R}^{k \times k},
\end{equation}
with $Q^TQ=I_k$ $R$ upper triangular. We have that:
\begin{align}
  \Phi^{\dagger} &= R^{-1} Q^T, \\
  \Phi^T \Phi &= R^T R.
\end{align}
We then write the first factor in \eqref{eq:tilde phi expl} as
\begin{align}
  (\Phi^T \Phi + \Delta)^{-1} 
  &= (R^T R + \Delta)^{-1} 
  \nonumber \\
  &= \left[
    R^T \left(
      I + R^{-T} \Delta R^{-1} 
    \right) R
  \right]^{-1}
  \nonumber \\
  &= R^{-1} \left[
    I + \sum_{l=1}^{\infty} (-1)^l
    \left( R^{-T} \Delta R^{-1} \right)^l
  \right]  R^{-T} 
  \nonumber \\
  &= (R^T R)^{-1} + S_{\Phi}
  \nonumber \\
  &= (\Phi^T \Phi)^{-1} + S_{\Phi},
  \label{eq:with Neumann}
\end{align}
where 
\begin{equation}
  S_{\Phi} = 
  R^{-1} \left[
    \sum_{l=1}^{\infty} (-1)^l
    \left( R^{-T} \Delta R^{-1} \right)^l
  \right]  R^{-T}
  \in \mathbb{R}^{k \times k},
    \label{eq:Amatrix}
\end{equation}
and in the second inequality in \eqref{eq:with Neumann}
we applied the Neumann series expansion to the 
matrix $I-R^{-T} \Delta R^{-1}$, 
which converges if 
\begin{equation}
    \|
        -R^{-T} \Delta R^{-1}       
    \|_2 < 1.
    \label{eq:cond Neumann}
\end{equation}
We will return to condition \eqref{eq:cond Neumann} 
at the end of this section.
We now substitute \eqref{eq:with Neumann} 
in \eqref{eq:tilde phi expl}, giving
\begin{align*}
  \tilde{\Phi}^{\dagger}
  &= \left[
    (\Phi^T \Phi)^{-1} + S_{\Phi}
  \right](\Phi^T + E^T)
  \nonumber \\
  &= \Phi^{\dagger} + (\Phi^T \Phi)^{-1} E^T 
    + S_{\Phi} \Phi^T + S_{\Phi} E^T,
\end{align*}
so we have that
\begin{equation}
  F^T = (\Phi^T \Phi)^{-1} E^T + S_{\Phi} \Phi^T + S_{\Phi} E^T,
  \label{eq:F expl}
\end{equation}
which is indeed $O(E)$, since $S_{\Phi} = O(\Delta)$ 
and $\Delta=O(E)$.
We next upper bound $\|S_{\Phi}\|_2$.  Firstly, note that, because
$\Phi^{\dagger} = \left(QR^{-1}\right)^T$, we have that:\footnote{
  To see the second equality in \eqref{eq:orthog}, for a
  matrix $Q \in \mathbb{R}^{m \times k}$ with $Q^TQ=I$ 
  and any matrix $A \in \mathbb{R}^{k \times k}$ we have that
  \begin{equation*}
    \|QA\|_2 = \sup_{\|v\|_2=1} \|QAv\|_2 
    = \sup_{\|v\|_2=1} \|Av\|_2 
    = \|A\|_2,
  \end{equation*}
  since 
  \begin{equation}
    \|QAv\|_2^2 = v^T A^T Q^T Q A v 
    = v^T A^T A v = \|Av\|_2^2.
  \end{equation}
}
\begin{align}
  \| \Phi^{\dagger} \|_2 = \| QR^{-1} \|_2 
  = \|R^{-1}\|_2.
  \label{eq:orthog}
\end{align}
Then, by using \eqref{eq:orthog}, norm submultiplicativity
and triangle inequality, from \eqref{eq:Amatrix} we have
\begin{equation}
  \| S_{\Phi} \|_2 \leq
  \| \Phi^{\dagger} \|^2_2
  \sum_{l=1}^{\infty}
  \| \Phi^{\dagger}\|^{2l}_2 \| \Delta \|^l_2.
  \label{eq:A first upper bnd}
\end{equation}
Now let $D$ be an upper bound on $\|\Delta\|_2$, obtained by
applying the triangle inequality in \eqref{eq:delta def}, so that
\begin{equation}
  \| \Delta \|_2 \leq D
  = 2 \| E \|_2 \| \Phi\|_2 + \|E\|^2_2.
  \label{eq:def D}
\end{equation}
Then, from \eqref{eq:A first upper bnd} we have
\begin{align}
  \|S_{\Phi}\|_2 &\leq 
  \| \Phi^{\dagger} \|^2_2
  \sum_{l=1}^{\infty}
  \| \Phi^{\dagger}\|^{2l}_2 D^l
  \nonumber \\
  &= \| \Phi^{\dagger} \|^2_2 \left(
    \frac{1}{ 
      1 - D \| \Phi^{\dagger}\|^2_2
    } -1
  \right)
  = \frac{
    D \| \Phi^{\dagger}\|^4_2
  }{
    1 - D \| \Phi^{\dagger}\|^2_2
  },
  \label{eq:bound norm A}
\end{align}  
where the series converges if $D\|\Phi^{\dagger}\|_2^2 < 1$,
in which case the denominator in the last fraction above is 
positive. We return to this condition at the end of 
the section.
We also know that\footnote{
  Using the SVD $\Phi = U \Sigma V^T$, we have
  $  
    \Phi^{\dagger} = (\Phi^T \Phi)^{-1} \Phi^T
    = (V \Sigma^2 V^T)^{-1} V \Sigma U^T 
    = V \Sigma^{-1} U^T
  $,
  so the conclusion follows.
}
\begin{equation}
  \| \Phi^{\dagger} \|_2 = \frac{1}{\sigma_{\min}(\Phi)}.
  \label{eq:norm phi dagger}
\end{equation}
By applying triangle inequality in \eqref{eq:F expl} and 
then using \eqref{eq:bound norm A} and the fact 
that $\|(\Phi^T\Phi)^{-1}\|_2 = 1/\sigma_{\min}^2(\Phi)
=\|\Phi^{\dagger}\|_2^2$ (from \eqref{eq:norm phi dagger}),
we obtain
\begin{equation}
  \|F\|_2 \leq 
  \|E\|_2  \|\Phi^{\dagger}\|_2^2
  + \frac{
    D \| \Phi^{\dagger}\|^4_2
  }{
    1 - D \| \Phi^{\dagger}\|^2_2
  } 
  \left( \| \Phi \|_2 + \|E\|_2 \right),
  \label{eq:norm F}
\end{equation}
where $D$ is given in \eqref{eq:def D}. 
It remains to establish an upper bound 
on $\|E\|_F$, and consequently on $\|E\|_2$.
The following lemma,
proved in Section~\ref{sec:proof Lemma bnd E_F}
gives us such a bound.

\begin{lemma}
    \textbf{(Upper bound on $\|E\|_F)$}
    Let $E = \tilde{\Phi} - \Phi$ for $\Phi$ and $\tilde{\Phi}$ 
    as defined in \eqref{eq:def Phi} and \eqref{eq:phi tilde def} 
    respectively for $t_j, \tilde{t}_j \in [0,1]$ 
    for $j=1,\ldots,k$. Then:
    \begin{equation}
        \|E\|_F 
        \leq
        \frac{4  
          e^{\frac{4}{\sigma^2}\max_j|\tilde{t}_j-t_j|}
          \sqrt{m}
        }{\sigma^2} 
          \|\mathbf{\tilde{t}}-\mathbf{t^*}\|_2.
        \label{eq:norm E}
    \end{equation}
    \label{lem:bound of E_F}
\end{lemma}

By using triangle inequality and 
norm submultiplicativity 
in \eqref{eq:a tilde}, and then 
substituting \eqref{eq:norm F} and \eqref{eq:norm E},
we obtain
\begin{align}
  \|\mathbf{a^*} - \mathbf{\tilde{a}}\|_2 
  &\leq 
    \|E\|_2 \|\Phi^{\dagger}\|_2 \|\mathbf{a^*}\|_2 
    + \|E\|_2^2 \|\Phi^{\dagger}\|_2^2 \|\mathbf{a^*}\|_2 
    \nonumber \\
    &\qquad + \frac{
      \|E\|_2 D \|\Phi^{\dagger}\|_2^4 
    }{
      1 - D \|\Phi^{\dagger}\|_2^2
    }
    (\|\Phi\|_2 + \|E\|_2)
    \|\mathbf{a^*}\|_2
  \nonumber \\
  & \leq  
  \frac{
    4 e^{\frac{4}{\sigma^2}\max_j|\tilde{t}_j-t_j|}
    \sqrt{m} \|\mathbf{a^*}\|_2
  }{
    \sigma^2 \sigma_{\min}(\Phi)
  } \|\mathbf{\tilde{t}}-\mathbf{t^*}\|_2
  + O(\|\mathbf{\tilde{t}}-\mathbf{t^*}\|_2^2),
  \label{eq:a bound expansion}
\end{align}
which is the bound given in Theorem~\ref{thm:a dep t}.
Note that because 
$\|E\|_2 = O(\|\tilde{\mathbf{t}} -\mathbf{t^*}\|_2)$
(see \eqref{eq:norm E}), 
the first
term is the only term that is 
$O(\|\tilde{\mathbf{t}} -\mathbf{t^*}\|_2)$
in the first inequality above, so
the other terms are included in the 
$O(\|\tilde{\mathbf{t}} -\mathbf{t^*}\|_2^2)$ 
term at the end.\footnote{
  In these terms, note that $\max_j |\tilde{t}_j-t_j| \leq 1$
  and therefore the notation 
  $O(\|\tilde{\mathbf{t}} -\mathbf{t^*}\|_2^2)$ 
  is correct.
}

Lastly, we return to condition \eqref{eq:cond Neumann}, 
which must be satisfied in order for the 
bound above to hold. By using norm submultiplicativity
and the bound on $\|\Delta\|_2$ from \eqref{eq:def D}, 
we obtain
\begin{equation}
    \| 
        {\Phi^{\dagger}}^T \Delta \Phi^{\dagger} 
    \|_2
    \leq 
    \| \Phi^{\dagger}\|^2_2 \|E\|^2_2
    + 2 \|\Phi\|_2 \|\Phi^{\dagger}\|^2_2 \|E\|_2
    \label{eq:rhs cond}
\end{equation}
and by requiring that the right hand side above is
less than one, we obtain a quadratic constraint on $\|E\|_2$, 
satisfied if
\begin{equation*}
    \|E\|_2 <
    \sigma_{\max}(\Phi)\left(
        \sqrt{
            1 + \frac{
                \sigma^2_{\min}(\Phi)
            }{
                \sigma^2_{\max}(\Phi) 
            }
        } -1
    \right).
\end{equation*}
By using the bound on $\|E\|_2$ from \eqref{eq:norm E}
with $\max_j |\tilde{t}_j-t_j| \leq 1$,
the above holds if 
\begin{equation}
  \delta_t
  <
  \frac{
      \sigma^2 \sigma_{\max}(\Phi)
  }{
      4e^{4/\sigma^2} \sqrt{m}
  }
  \left(
      \sqrt{
          1 + \frac{
              \sigma_{\min}^2(\Phi)
          }{
              \sigma_{\max}^2(\Phi)
          } 
      } -1
  \right),
  \nonumber
\end{equation} 
which is the condition \eqref{eq:cond t} in the statement
of the theorem.
Note that by imposing this, we also ensure that the condition
for the series in \eqref{eq:bound norm A} to converge holds, 
since $D\|\Phi^{\dagger}\|_2^2$ is equal to the right hand side
of \eqref{eq:rhs cond}.

\subsubsection{Proof of Lemma~\ref{lem:bound of E_F} 
  (Bound of $\|E\|_F$)}
  \label{sec:proof Lemma bnd E_F}

Since $E = \tilde{\Phi} - \Phi$, 
for $\tilde{t}_j$ being a perturbation of $t_j$, we have that
\begin{equation*}
  |E_{ij}| 
  = \left| 
    e^{-\frac{(s_i-\tilde{t}_j)^2}{\sigma^2}} 
    - e^{-\frac{(s_i-t_j)^2}{\sigma^2}} 
  \right| 
  = e^{-\frac{(s_i-t_j)^2}{\sigma^2}} \left|
    e^{\frac{1}{\sigma^2}
      \left[(s_i-t_j)^2 - (s_i-\tilde{t}_j)^2 \right] } - 1
  \right|.
\end{equation*}
Then the exponent can be written as
\begin{equation*}
  \left| \frac{1}{\sigma^2}
  \left[
    (s_i-t_j)^2 - (s_i-\tilde{t}_j)^2 
  \right] \right|
  = 
  \left| \frac{1}{\sigma^2}
  \left[
    2s_i(\tilde{t}_j-t_j) + (t_j+\tilde{t}_j)(t_j-\tilde{t}_j)
  \right] \right|
  \leq  
    \frac{4|\tilde{t}_j - t_j|}{\sigma^2},
\end{equation*}
where we used that $ s_i, \tilde{t}_j, t_j \in [0,1]$, so
\begin{equation*}
  e^{-\frac{4}{\sigma^2} |\tilde{t}_j - t_j|}
  \leq
  e^{\frac{1}{\sigma^2}
    \left[(s_i-t_j)^2 - (s_i-\tilde{t}_j)^2 \right] } 
  \leq
  e^{\frac{4}{\sigma^2}{|\tilde{t}_j - t_j|}},
\end{equation*}
which implies that
\begin{align*}
    |E_{ij}| 
    &\leq 
  \left|
    e^{\frac{1}{\sigma^2}
      \left[(s_i-t_j)^2 - (s_i-\tilde{t}_j)^2 \right] } 
    - 1
  \right|
  \nonumber \\
  &\leq 
  \max\left\{
    1 - e^{-\frac{4}{\sigma^2}|\tilde{t}_j-t_j|}, 
    e^{\frac{4}{\sigma^2}|\tilde{t}_j - t_j|} -1
  \right\}
  \nonumber \\
  &= e^{\frac{4}{\sigma^2}|\tilde{t}_j - t_j|} -1
  \nonumber \\
  &= \frac{4}{\sigma^2} | \tilde{t}_j - t_j| 
  \cdot e^{\xi},
\end{align*} 
for some $\xi \in 
\left[
  -\frac{4}{\sigma^2} | \tilde{t}_j - t_j|, 
  \frac{4}{\sigma^2} | \tilde{t}_j - t_j|
\right]$ 
and where in the first inequality we have used that 
$e^{-\frac{(s_i-t_j)^2}{\sigma^2}} \leq 1$. 
Then
\begin{align*}
  |E_{ij}| &\leq
  \frac{4}{\sigma^2} | \tilde{t}_j - t_j| 
  \cdot e^{
    \frac{4}{\sigma^2} | \tilde{t}_j - t_j|
  }, 
  %
\end{align*}
and we conclude that
\begin{align}
  \|E\|_2 \leq 
  \|E\|_F 
  &\leq
  \sqrt{
    \sum_{i=1}^m \sum_{j=1}^k
    \left(
      \frac{
        4 e^{\frac{4}{\sigma^2}|\tilde{t}_j-t_j|  } 
      }{\sigma^2}
    \right)^2
    | \tilde{t}_j - t_j|^2
  }
  \nonumber \\
  &\leq
  \frac{4  
    e^{\frac{4}{\sigma^2}\max_j|\tilde{t}_j-t_j|}
    \sqrt{m}
  }{\sigma^2} 
    \|\mathbf{\tilde{t}}-\mathbf{t^*}\|_2,
\end{align}
provided that $\tilde{t}_j,t_j \in [0,1]$
for all $j=1,\ldots,k$.

\subsection{Proof of Theorem~\ref{thm:pert noisy}
  (Dependence of $\|\lambda-\lambda^*\|_2$ on the noise $w$)}
\label{sec:pert noisy proof}

We apply the quantitative implicit function theorem
to the function $\bar{F}$ defined in \eqref{eq:F bar}. 
First, note that in the bound \eqref{eq:pert noisy}, which
we want to prove, we only need to consider the $k$ non-fixed
entries of $\lambda$ and $\lambda^*$, as the error in the other
entries zero. Therefore, in this proof we will
only work with the vectors of non-fixed entries
$\bar{\lambda}, \bar{\lambda}^* \in \mathbb{R}^k$, but we will abuse
the notation for simplicity and write $\lambda$ and $\lambda^*$ 
respectively. Similarly, we will write $w$ to denote
the vector $w_{\theta} \in \mathbb{R}^{2k}$ 
corresponding to the $2k$ entries
of the noise vector $w \in \mathbb{R}^m$
and $s_1,\ldots,s_{2k}$ to denote
$s_{\theta_{1}},\ldots,s_{\theta_{2k}}$.
The partial derivatives of $\bar{F}$ 
from Definition~\ref{def:F bar} are:
\begin{align}
  \partial_{\lambda_l} \bar{F}_j 
  &=
  - \sum_{i=1}^k \nu_i \phi'(t_i(\lambda)-s_j) 
  \partial_{\lambda_l} t_i(\lambda)
  \quad l=1,\ldots,k,
  \quad j=1,\ldots,2k,
  \label{eq:dF d lam}
  \\
  \partial_{\nu_l} \bar{F}_j
  &= -\phi(t_l(\lambda) -s_j),
  \quad l=1,\ldots,k,
  \quad j=1,\ldots,2k,
  \label{eq:dF d nu}
  \\
  \partial_{w_l} \bar{F}_j
  &= \begin{cases}
    1, \quad \text{if} \quad l=j,
    \\
    0, \quad \text{otherwise},
  \end{cases}
  \quad l,j=1,\ldots,2k.
  \label{eq:dF d eps}
\end{align}

Let $\gamma = [\lambda,\nu]^T$ 
and $\gamma^* = [\lambda^*,a]^T$, so that we can 
write $\bar{F}([\lambda,\nu]^T,w)$ 
as $\bar{F}(\gamma,w)$ and $\bar{F}(\gamma^*,0) = 0$.
In order to apply the implicit function theorem, 
the following conditions must be satisfied:
\begin{enumerate}
  \item $\partial_{\gamma} \bar{F}(\gamma^*,0)$ is invertible,

  \item We choose the radius $\delta_{\gamma}$ 
    of the ball $V_{\delta_\gamma}$ around $\gamma$ where the result of the 
    quantitative implicit function theorem holds:
    \begin{equation}
      \sup_{(\gamma,w) \in V_{\delta_\gamma}}
      \left\|
        I - \left[
          \partial_{\gamma} \bar{F}(\gamma^*,0)
        \right]^{-1}
        \partial_{\gamma} \bar{F}(\gamma,w)
      \right\|_2 \leq \frac12,
    \end{equation}

  \item The radius $\delta_{w}$ of the ball 
    around $w^*=0$ that contains $w$ is:
    \begin{equation}
      \delta_{w} = (2M_{w}B_{\delta_{\gamma}})^{-1} \delta_{\gamma},
    \end{equation}
    where
    \begin{align}
      B_{\delta_{\gamma}} &= \sup_{(\gamma,w)\in V_{\delta_{\gamma}}}
      \|
        \partial_{w} \bar{F}(\gamma,w)
      \|_2,
      \\
      M_{w} &= \|
        \partial_{\gamma} \bar{F}(\gamma^*,0)^{-1}
      \|_2.
    \end{align}
\end{enumerate}
The first condition is also one of the conditions in the theorem,
and it has been discussed in Section~\ref{sec:discussion}. We now need
to establish the two radii for the balls of the perturbations.

\subsubsection*{Perturbation radii}

Before proceeding to calculating the radii of the balls where the 
implicit function theorem holds, we state the following lemma,
which allows us to write the Jacobian of $\bar{F}$ with respect to the first
variable as a sum of the Jacobian evaluated 
at $(\gamma^*,w^*)=([\lambda^*,a]^T,0)$ and a perturbation matrix, 
whose norm is bounded explicitly.
The proof of Lemma~\ref{lem:jacobian pert} is given
in Section~\ref{sec:proof pert jac}.

\begin{lemma}
  \label{lem:jacobian pert}
  \textbf{(Bound on the perturbation of the Jacobian of $\bar{F}$)}
  Let $J(\lambda,\nu,w)$ be the Jacobian of $\bar{F}(\gamma,w)$
  with respect to $\gamma = [\lambda,\nu]^T$
  and $\delta_{\gamma}$ an upper bound on the perturbation 
  of $\gamma^*=[\lambda^*,a]^T$, namely:
  \begin{equation*}
    \left\|
      \begin{bmatrix}
        \lambda - \lambda^* \\
        \nu - a
      \end{bmatrix}
    \right\|_2 \leq \delta_{\gamma}.
  \end{equation*}
  Then:
  \begin{equation}
    J(\lambda,\nu,w) = J(\lambda^*,a,0) + E,
  \end{equation}
  with 
  \begin{equation}
    \|E\|_F \leq  P(k,\sigma,\Pi,\tau,C_{t^*}) \cdot \delta_{\gamma},
  \end{equation}
  where:
  \begin{align}
    P(k,\sigma,\Pi,&\tau,C_{t^*}) 
      =\sqrt{2} k
      \Bigg[
        \frac{1}{\sigma^2} 
        \left(
          2 \sqrt{k} C_{t^*}^2 \Pi
          + 4 k C_{t^*} \bar{\Delta}_2 \tau \Pi
          + 2 C_{t^*}
        \right)
      \nonumber \\
      &+ \frac{1}{\sigma} 
      \left(
        \frac{\sqrt{2k} C_{t^*}}{\sqrt{e}}
        + 4 \sqrt{k} C_{t^*}^2 \Pi
        + \frac{2\sqrt{2} \bar{\Delta}_2 \Pi}{\sqrt{e}}
        + 8 k C_{t^*} \bar{\Delta}_2 \tau \Pi
        + \frac{\sqrt{2k} \bar{\Delta}_2 \Pi}{\sqrt{e}}
      \right)
      \Bigg],
    \label{eq:E_F bnd}
  \end{align}
  for $\|\lambda - \lambda^*\| \leq \delta_{\lambda}$,
  where $\delta_{\lambda}$ and $C_{t^*}$ are
  given in \eqref{eq:delta_lambda} and \eqref{eq:C_t}
  respectively in Theorem~\ref{thm:t dep lambda}
  and $\bar{\Delta}_2$ is given in \eqref{eq:D2 bnd}
  in the proof.
\end{lemma}

We can now use Lemma~\ref{lem:jacobian pert} to write
\begin{equation}
  \partial_{\gamma} \bar{F}(\gamma,w) = \partial_{\gamma} \bar{F}(\gamma^*,0)+ E,
\end{equation}
then
\begin{align}
  I - \left[
    \partial_{\gamma} \bar{F}(\gamma^*,0)
  \right]^{-1}
  \partial_{\gamma} \bar{F}(\gamma,w)
  &=
  I - \left[
    \partial_{\gamma} \bar{F}(\gamma^*,0)
  \right]^{-1}
  \left[
    \partial_{\gamma} \bar{F}(\gamma^*,0) + E
  \right]
  \nonumber \\
  &= 
  -\left[
    \partial_{\gamma} \bar{F}(\gamma^*,0)
  \right]^{-1} E,
\end{align}
so 
\begin{align}
  \left\|
    I - \left[
      \partial_{\gamma} \bar{F}(\gamma^*,0)
    \right]^{-1}
    \partial_{\gamma} \bar{F}(\gamma,w)
  \right\|_2 
  &\leq
  \left\|
    \left[
      \partial_{\gamma} \bar{F}(\gamma^*,0)
    \right]^{-1}
  \right\|_2 \cdot
  \| E \|_F
  \nonumber \\ &\leq
  \frac{
    \|E\|_F
  }{
    \sigma_{\min}(\partial_{\gamma} \bar{F}(\gamma^*,0))
  }
  \nonumber \\ &\leq
  \frac{
    P(k,\sigma,\Pi,\tau,C_{t^*}) 
  }{
    \sigma_{\min}(\partial_{\gamma} \bar{F}(\gamma^*,0))
  } \cdot \delta_{\gamma},
  \label{eq:condition delta x noise}
\end{align}
where $P \cdot \delta_{\gamma}$ is the upper bound on $\|E\|_F$
given in \eqref{eq:E_F bnd}.

Therefore, from the condition that the right-hand side of
the last inequality is less than or equal to $\frac12$,
we choose the radius $\delta_{\gamma}$ to be:
\begin{equation}
  \delta_{\gamma} =\frac{
    \sigma_{\min}(\partial_{\gamma} \bar{F}(\gamma^*,0))
  }{
    2 P(k,\sigma,\Pi,\tau,C_{t^*}) 
  }.
  \label{eq:delta x}
\end{equation}
Using \eqref{eq:dF d eps}, we have that
\begin{equation}
  B_{\delta_{\gamma}} = 1.
\end{equation}
Then 
\begin{equation}
  M_{w} = \frac{1}{\sigma_{\min}(\partial_{\gamma} \bar{F}(\gamma^*,0) )},
\end{equation}
so, using \eqref{eq:delta x}, we obtain
\begin{equation}
  \delta_{w} = \frac{
    \sigma_{\min}(\partial_{\gamma} \bar{F}(\gamma^*,0))^2
  }{
    4 P(k,\sigma,\Pi,\tau,C_{t^*}) 
  }.
\end{equation}

\subsubsection*{Applying the quantitative implicit function theorem}

Having calculated the radii where the quantitative implicit 
function theorem holds, we apply it to obtain:
\begin{equation}
  \partial_{w} g(w) = 
  - \left[
      \partial_1 \bar{F}(g(w),w)
    \right]^{-1},
  \label{eq:partial g eps}
\end{equation}
where $\partial_1$ is the partial derivative with
respect with the first argument and $g(w)$
gives the dependence of $[\lambda,\nu]^T$ on $w$.
Specifically, we write:
\begin{align}
  &\lambda_i(w) = g_i(w)
  \quad \text{ for } \quad
  i = 1,\ldots,k,
  \label{eq:lam_g}
  \\
  &\nu_i(w) = g_{k+i}(w)
  \quad \text{ for } \quad
  i = 1,\ldots,k.
  \label{eq:nu_g}
\end{align}
Let 
$
  J(\lambda,\nu,w) = 
  \partial_1 \bar{F}([\lambda,\nu]^T,w)
$, 
where $\lambda = \lambda(w)$
and $\nu = \nu(w)$ 
by \eqref{eq:lam_g} and \eqref{eq:nu_g}.
Lemma~\ref{lem:jacobian pert} gives
\begin{equation}
  J(\lambda,\nu,w) = 
  J(\lambda^*,a,0) + E,
\end{equation} 
so $E$ is the perturbation
of $J(\lambda^*,a,0)$ due to perturbed
$\lambda, \nu, w$ and a bound
on $\|E\|_F$ is given in the lemma.

We will now use the following result,
proved in Section~\ref{sec:sig pert proof},
which enables us to make use of the upper bound 
on the norm of the perturbation given
by Lemma~\ref{lem:jacobian pert} in
order to lower bound the smallest singular 
value of~$J$.

\begin{lemma}
  \label{lem:sig pert}
  Let $J \in \mathbb{R}^{m \times n}$. 
  If $J = A+E$, then 
  \begin{equation*}
    \sigma_{\min}(J) \geq \sigma_{\min}(A) - \|E\|_F.
  \end{equation*}
\end{lemma}

By applying Lemma~\ref{lem:sig pert}, we have that:
\begin{align}
  \frac{1}{\sigma_{\min}(J(\lambda,\nu,w))} &\leq 
  \frac{1}{\sigma_{\min}(J(\lambda^*,a,0)) - \|E\|_F}
  = \frac{1}{
    \sigma_{\min}(J(\lambda^*,a,0)) \left( 1 -
      \frac{\|E\|_F}{\sigma_{\min}(J(\lambda^*,a,0))}
    \right)
  }
  \nonumber \\ &\leq
  \frac{1}{\sigma_{\min}(J(\lambda^*,a,0))} \cdot
  \left(
    1 + \frac{2\|E\|_F}{\sigma_{\min}(J(\lambda^*,a,0))}
  \right),
  \label{eq:sig A E}
\end{align}
for 
\begin{equation}
  \frac{\|E\|_F}{\sigma_{\min}(J(\lambda^*,a,0))} \leq \frac12,
\end{equation}
where we used the fact that
$ (1-x)^{-1} \leq 1 + 2x$ for $x \in [0,\frac12]$.
Note that the condition above is the same as the 
condition that the right hand side 
of \eqref{eq:condition delta x noise} is less than 
or equal to $\frac12$, which
is satisfied for our choice of $\delta_{\gamma}$ and $\delta_{w}$.

From \eqref{eq:partial g eps}
and \eqref{eq:sig A E}, we have that:
\begin{align}
  \| \partial_{w} g(w)\|_2 
  &= \frac{1}{\sigma_{\min}(J(\lambda^*,a,0) + E)}
  \nonumber \\ &\leq
    \frac{1}{\sigma_{\min}(J(\lambda^*,a,0))}
    \left(
      1 + \frac{2}{\sigma_{\min}(J(\lambda^*,a,0))} 
      \cdot \|E\|_F
    \right)
  \nonumber \\ &\leq
  \frac{2}{
    \sigma_{\min}(J(\lambda^*,a,0))
  },
  \label{eq:Hessian bound}
\end{align}
where $\|E\|_F$ is upper bounded in \eqref{eq:E_F bnd} and
$w$, $\lambda$ and $\nu$ satisfy
\begin{equation*}
  \| w \|_2 \leq \delta_{w},
  \quad
  \| \lambda - \lambda^* \|_2 \leq~ \delta_{\gamma},
  \quad
  \| \nu - a\|_2 \leq \delta_{\gamma}.
\end{equation*}
The first-order Taylor expansion of $g(w)$ 
around $w=0$ is:
\begin{equation}
  g(w) = g(0) + \partial_{w} g(w_{\delta})^T w,
  \label{eq:g taylor}
\end{equation}
for some $w_{\delta}$ on the segment between the 
zero vector and $w$. Noting that $g(w)$ 
is our notation for the vector:
\begin{equation}
  g(w) = \begin{bmatrix}
    \lambda(w) \\
    \nu(w)
  \end{bmatrix},
\end{equation}
with $\lambda(0) = \lambda^*$ and $\nu(0) = a$, 
from \eqref{eq:g taylor} we have that:
\begin{align}
  \left\| \begin{bmatrix}
    \lambda(w) - \lambda^* \\
    \nu(w) - a
  \end{bmatrix} \right\|_2
  &= 
  \left\|
    \partial_{w}g(w_{\delta})^T
    w 
  \right\|_2
  \nonumber \\
  &\leq \|\partial_{w}g(w_{\delta})\|_2 
  \cdot \|w\|_2, 
\end{align}
for $w$, $\lambda$ and $\nu$ such that
\begin{equation*}
  \| w \|_2 \leq \delta_{w},
  \quad
  \| \lambda - \lambda^* \|_2 \leq \delta_{\gamma},
  \quad
  \| \nu - a \|_2 \leq \delta_{\gamma},
\end{equation*}
where we use the bound from \eqref{eq:Hessian bound}.

\subsubsection{Proof of Lemma~\ref{lem:sig pert}}
  \label{sec:sig pert proof}

We have that 
\begin{align*}
  \sigma(J) &= 
  \min_{\|v\|_2=1} \max_{\|u\|_2=1} u^T(A+E) v
  \\ &\geq
  \min_{\|v\|_2=1} \max_{\|u\|_2=1} u^T A v 
  - \max_{\|v\|_2=1} \max_{\|u\|_2=1} u^T E v 
  \\ &\geq
  \sigma_{\min}(A) - \|E\|_F.
\end{align*}

\subsubsection{Proof of Lemma~\ref{lem:jacobian pert}
  (Bound on the perturbation of the Jacobian of $\bar{F}$)}
  \label{sec:proof pert jac}

Let 
$
  J(\lambda,\nu,w) = 
  \partial_1 \bar{F}([\lambda,\nu]^T,w)
$, 
where $\lambda = \lambda(w)$ 
and $\nu = \nu(w)$ 
by \eqref{eq:lam_g} and \eqref{eq:nu_g},
and we want to write $J$ in the form
\begin{equation}
  J(\lambda,\nu,w) = 
  J(\lambda^*,a,0) + E
  \label{eq:J J E}
\end{equation} 
i.e. $E$ is the perturbation
of $J(\lambda^*,a,0)$ due to perturbed
$\lambda, \nu, w$.
In order to apply Lemma~\ref{lem:sig pert}, we need an
upper bound on $\|E\|_F$, so we need to 
upper bound each entry of $E$. Let 
\begin{equation}
  J = \left[
    J_1 J_2
  \right],
\end{equation}
where $J_1$ corresponds to the 
terms \eqref{eq:dF d lam} and $J_2$ 
to the terms \eqref{eq:dF d nu}
and 
\begin{equation}
  E = \left[
    E_1  E_2
  \right]
\end{equation}
the corresponding perturbation terms.

\subsection*{Entries in $J_1$}

For $i=1,\ldots,k$ and $j=1,\ldots,2k$:
\begin{align}
  J_{1_{j,i}} &= 
  - \sum_{p=1}^k (\nu_p -a_p+a_p) 
    \phi'(t_p(\lambda) -t_p^* + t_p^*   - s_j) 
    \partial_{\lambda_i} t_p(\lambda)
  \nonumber \\ &=
  - \sum_{p=1}^k 
    \partial_{\lambda_i} t_p(\lambda)
    \left[
      a_p \phi'(t_p^* -s_j + t_p(\lambda)-t_p^*)
      + (\nu_p-a_p) \phi'(t_p^* - s_j + t_p(\lambda)-t_p^*)
    \right]
  \nonumber \\ &=
  - \sum_{p=1}^k 
    \partial_{\lambda_i} t_p(\lambda)
    \big[
      a_p \phi'(t_p^* -s_j) 
      + a_p(t_p(\lambda)-t_p^*)\phi''(\xi_{j,p})
      + (\nu_p-a_p) \phi'(t_p^* - s_j) 
      \nonumber \\
      &\quad+ (\nu_p-a_p)(t_p(\lambda)-t_p^*)\phi''(\xi_{j,p})
    \big]
  \nonumber \\ &=
  - \sum_{p=1}^k 
    \partial_{\lambda_i} t_p(\lambda)
    \left(
      a_p \phi'(t_p^* -s_j) + \Delta_{1_{j,p}}
    \right),
    \label{eq:J1ij}
\end{align}
where
\begin{equation}
  \Delta_{1_{j,p}} =
    a_p(t_p(\lambda)-t_p^*)\phi''(\xi_{j,p})
    + (\nu_p-a_p) \phi'(t_p^* - s_j) 
    + (\nu_p-a_p)(t_p(\lambda)-t_p^*)\phi''(\xi_{j,p}),
\end{equation}
for some $\xi_{j,p} \in [t_p^*  -s_j -|t_p-t_p^*|,t_p^*-s_j+|t_p-t_p^*|]$.
The factor involving the partial derivative 
in \eqref{eq:J1ij} has the same form as \eqref{eq:qift res}
so in order to bound it 
we write the Taylor expansion 
of \eqref{eq:qift res} around $\lambda^*$:
\begin{equation}
  \partial_{\lambda} t_p(\lambda)
  = \partial_{\lambda} t_p(\lambda^*)
  + \partial^2_{\lambda\lambda} t_p (\lambda_{\delta})
    \left(\lambda - \lambda^*\right),
  \label{eq:dtp d lam exp}
\end{equation}
for some $\lambda_{\delta}$ on the segment 
between $\lambda$ and $\lambda^*$. 
By using \eqref{eq:qift res} with \eqref{eq:partial t}
and \eqref{eq:partial lambda}, the entry $i,l$ 
in the Hessian matrix 
$H = \partial^2_{\lambda\lambda} t_p (\lambda_{\delta})$ is
\begin{equation}
  \left(H^T\right)_{i,l}
  = \frac{ F_{i,l}(\lambda) }{
    \left[
      \sum_{j=1}^m \lambda_j
      \phi''(t_p(\lambda)-s_j)
    \right]^2
  },
  \label{eq:H entry}
\end{equation}
for $i,l = 1,\ldots,k$, where
\begin{align}
  F_{i,l}(\lambda) = 
  &- \phi''(t_p(\lambda)-s_i) \partial_{\lambda_l} t_p(\lambda) 
  \sum_{j=1}^m \lambda_j \phi''(t_p(\lambda)-s_j)
  \nonumber \\
  & + \phi'(t_p(\lambda)-s_i) 
  \left(
    \sum_{j=1}^m \lambda_j \phi'''(t_p(\lambda)-s_j) 
    \partial_{\lambda_l} t_p(\lambda)
    + \phi''(t_p(\lambda)-s_l)
  \right).
  \label{eq:Fil def}
\end{align}
Note that in the denominator \eqref{eq:H entry} 
we use all $m$ entries of $\lambda$ and samples
due to how we defined the function from \eqref{eq:qift res},
and the same is true for the sums in \eqref{eq:Fil def}.
From \eqref{eq:dtp d lam exp} and \eqref{eq:Fil def}, we then write:
\begin{equation}
  \partial_{\lambda_i} t_p(\lambda) 
  = \partial_{\lambda_i} t_p(\lambda^*)
  + \Delta_{2_{i,p}},
\end{equation}
where
\begin{equation}
  \Delta_{2_{i,p}} = 
  \sum_{l=1}^{k} \frac{
    (\lambda_l - \lambda^*_l)
    F_{i,l}(\lambda_{\delta})
  }{
    \left[
      \sum_{j=1}^m \lambda_j
      \phi''(t_p(\lambda_{\delta})-s_j)
    \right]^2
  }.
\end{equation}
Note that $l$ goes up to $k$ because we only 
work with $k$ entries in $\lambda$.
Therefore, we have that:
\begin{equation}
  J_{1_{j,i}} = 
  - \sum_{p=1}^k 
  \left(
    \partial_{\lambda_i} t_p(\lambda^*)
    + \Delta_{2_{i,p}}
  \right)
  \left(
      a_p \phi'(t_p^* -s_j) + \Delta_{1_{j,p}}
    \right),
  \label{eq:J1 factors}
\end{equation}
where 
\begin{align}
  \Delta_{1_{j,p}} &= O(|t_p-t_p^*| + |\nu_p-a_p|), \\
  \Delta_{2_{i,p}} &= O(\|\lambda-\lambda^*\|_2),
\end{align}
for $i=1,\ldots,k$, $j=1,\ldots,2k$ and $p=1,\ldots,k$.
The next step now is to upper bound 
$|\Delta_{1_{j,p}}|$ and $|\Delta_{2_{i,p}}|$.

\subsubsection*{Bounding $\Delta_{1_{j,p}}$}

By the triangle inequality, we have that:
\begin{align}
  | \Delta_{1_{j,p}} | &\leq 
    |a_p| |t_p(\lambda)-t_p^*| |\phi''(\xi_{j,p})|
    + |\nu_p-a_p| | \phi'(t_p^* - s_j) |
    + |\nu_p-a_p| |t_p(\lambda)-t_p^*| |\phi''(\xi_{j,p})|
  \nonumber \\ &\leq
  |a_p| |t_p(\lambda)-t_p^*| \frac{2}{\sigma^2}
  + |\nu_p - a_p| \frac{\sqrt{2}}{\sqrt{e}\sigma}
  + |\nu_p - a_p| |t_p(\lambda) - t_p^*| \frac{2}{\sigma^2}
  =: \bar{\Delta}_{1_p},
  \label{eq:D1p bound}
\end{align}
for $j=1,\ldots,2k$ and $p=1,\ldots,k$, 
where we have used the maxima of the Gaussian and its 
derivatives given in footnote \ref{fn:derivs}.

\subsubsection*{Bounding $\Delta_{2_{i,p}}$}

By applying the Cauchy-Schwartz inequality, we have that:
\begin{align}
  |\Delta_{2_{i,p}}| &\leq
  \frac{1}{
    \left|
      \sum_{j=1}^m \lambda_j
      \phi''(t_p(\lambda_{\delta})-s_j)
    \right|^2
  } \cdot
  \left\| \begin{bmatrix}
    F_{i,1}(\lambda_{\delta}) \\
    \vdots \\
    F_{i,k}(\lambda_{\delta})
  \end{bmatrix} \right\|_2 \cdot
  \|\lambda - \lambda^* \|_2
  \label{eq:D2 bnd almost}
\end{align}
We now bound $|F_{i,l}|$ for $i,l=1,\ldots,k$:
\begin{align}
  |F_{i,l}(\lambda_{\delta})| &\leq
  |\phi''(t_p(\lambda_{\delta})-s_i)| 
  \left| 
    \partial_{\lambda_l} t_p(\lambda_{\delta})
  \right|
  \left|
    \sum_{j=1}^m \lambda_j \phi''(t_p(\lambda_{\delta})-s_j)
  \right|
  \nonumber \\
  &\qquad + |\phi'(t_p(\lambda_{\delta})-s_i) |
  \left(
    \left| 
      \partial_{\lambda_l} t_p(\lambda_{\delta})
    \right|
    \left|
      \sum_{j=1}^m \lambda_j \phi'''(t_p(\lambda_{\delta})-s_j)
    \right|
    + |\phi''(t_p(\lambda_{\delta})-s_l)|
  \right)
  \nonumber \\ &\leq
  \frac{2 C_{t^*}}{\sigma^2}
  \| \lambda_{\delta} \|_2 \cdot
  \left\|
    \left[ 
      \phi''(t_p(\lambda_{\delta})-s_j) 
    \right]_{j=1}^m
  \right\|_2 
  \nonumber \\
  &\qquad + \frac{\sqrt{2}}{\sqrt{e}\sigma}
  \left(
    C_{t^*}
    \|\lambda_{\delta}\|_2 \cdot
    \left\|
      \left[
        \phi'''(t_p(\lambda_{\delta})-s_j)
      \right]_{j=1}^m
    \right\|_2
    + \frac{2}{\sigma^2}
  \right)
  \nonumber \\ &\leq
  \frac{2 C_{t^*}}{\sigma^2}
  \| \lambda_{\delta} \|_2 \cdot
  \frac{2\sqrt{m}}{\sigma^2}
  + \frac{\sqrt{2}}{\sqrt{e}\sigma}
  \left(
    C_{t^*}
    \|\lambda_{\delta}\|_2 \cdot
    \frac{c\sqrt{m}}{\sigma^3}
    + \frac{2}{\sigma^2}
  \right),
\end{align}
where we used the Cauchy-Schwartz inequality,
the bounds in footnote \ref{fn:derivs} and $C_{t^*}$
from \eqref{eq:C_t}. Therefore, the above inequality
holds for $\lambda_{\delta} \in \mathcal{B}(\lambda^*,\delta_{\lambda})$
with $\delta_{\lambda}$ from \eqref{eq:delta_lambda}.

The final bound on $|F_{i,l}|$ is
\begin{equation}
  |F_{i,l}| \leq 
  \frac{c_2 C_{t^*} m \tau}{\sigma^4}
  + \frac{2\sqrt{2}}{\sqrt{e}\sigma^3},
  \label{eq:Fil bnd}
\end{equation}
where $c_2 = 4 + \frac{c\sqrt{2}}{\sqrt{e}} \approx 7.3484$,
for $i,l=1,\ldots,k$ and 
we used $\|\lambda_{\delta}\|_2 \leq \tau \sqrt{m}$.

The next step is to obtain a lower bound on the denominator
in \eqref{eq:D2 bnd almost}.
By adding and subtracting $\lambda_j^*$ to $\lambda_j$ and 
applying the reverse triangle inequality, we obtain:
\begin{align}
  \left|
    \sum_{j=1}^m \lambda_j \phi''(t_p(\lambda)-s_j)
  \right|
  &\geq
  \left|
    \sum_{j=1}^m \lambda_j^* 
    \phi''(t_p(\lambda)-s_j)
  \right|
  - \left|
    \sum_{j=1}^m (\lambda_j - \lambda_j^*) 
    \phi''(t_p(\lambda)-s_j)
  \right| \label{eq:bprime}
  \\ 
  &\geq
  |q''(t^*)|\left[ 
    1 - \frac{c \|\lambda^* \|_2}{4\sigma + 2c \|\lambda^* \|_2}    
  \right]
  - \frac{2 \sqrt{m} \|\lambda-\lambda^*\|_2}{\sigma^2},
  \label{eq:den temp}
\end{align}
where the first term on the right hand side on \eqref{eq:bprime}
has the same form as $B$ in \eqref{eq:B def}, and therefore
on the next line we use the bound in \eqref{eq:bound B}.
For the second term, we apply the Cauchy-Schwartz inequality
and the bound in footnote \ref{fn:derivs}, where the 
constant $c \approx 3.9036$ is obtained.  
The last inequality above holds under the condition that the
right hand side is positive.
We set the stronger condition that the right hand side
of~\eqref{eq:den temp} is greater than or equal to one:
\begin{equation}
  \left|
    \sum_{j=1}^m \lambda_j \phi''(t_p(\lambda)-s_j)
  \right|
  \geq 1,
  \label{eq:D2 den bnd}
\end{equation}
which is satisfied if: 
\begin{equation}
  |q''(t^*)| \geq \left( 1  
    + \frac{2 \sqrt{m} \|\lambda-\lambda^*\|_2}{\sigma^2}
  \right)
  \cdot
  \frac{4\sigma + 2 c \|\lambda^* \|_2}{4\sigma + c \|\lambda^* \|_2}.
  \label{eq:q cond temp}
\end{equation}
Then, using the box constraint $\|\lambda\|_{\infty} \leq \tau$ 
and the fact that $\|\lambda-\lambda^*\|_2 \leq 2\tau\sqrt{m}$
from the triangle inequality, and the fact that
$
\frac{4\sigma + 2 c \|\lambda^* \|_2}{4\sigma + c \|\lambda^* \|_2}
\leq 2
$
the condition~\eqref{eq:q cond temp} is satisfied if:
\begin{equation}
  |q''(t^*)| \geq 2 \left(
    1 + \frac{4 m \tau}{\sigma^2} 
  \right).
  \label{eq:cond q}
\end{equation}
By combining \eqref{eq:D2 bnd almost}, \eqref{eq:Fil bnd} 
and \eqref{eq:D2 den bnd}, we obtain
the final bound on $|\Delta_{2_{i,p}}|$:
\begin{equation}
  | \Delta_{2_{i,p}} | \leq
  \bar{\Delta}_2 \cdot
  \|\lambda - \lambda^*\|_2,
\end{equation}
for $i=1,\ldots,k$ and $p=1,\ldots,k$,
where
\begin{equation}
  \bar{\Delta}_2 =
  \frac{\sqrt{k}}{\sigma^4} \left(
    c_2 C_{t^*} m \tau
  + \frac{2\sqrt{2}}{\sqrt{e}} \sigma
  \right),
  \label{eq:D2 bnd}
\end{equation}
and $c \approx 3.9036$, $c_2 = 4 + \frac{c\sqrt{2}}{\sqrt{e}} \approx 7.3484$.

Therefore, from \eqref{eq:J1 factors} and by using
the definitions of $\bar{\Delta}_{1_p}$ 
and $\bar{\Delta}_2$ from \eqref{eq:D1p bound} 
and \eqref{eq:D2 bnd} respectively, we have that:
\begin{align}
  |E_{1_{j,i}}| &= \left|
    \sum_{p=1}^k 
    \partial_{\lambda_i} t_p(\lambda^*) \Delta_{1_{j,p}}
    + a_p \phi'(t_p^*-s_j) \Delta_{2_{i,p}}
    + \Delta_{1_{j,p}} \Delta_{2_{i,p}} \right|
    \nonumber \\ 
  &\leq
    C_t \sum_{p=1}^k \bar{\Delta}_{1_p}
    + \|\lambda - \lambda^*\|_2 \bar{\Delta}_2 
    \| a \|_2 \cdot
    \left\|
      \left[\phi'(t_p^*-s_j)\right]_{p=1}^k
    \right\|_2
    + \|\lambda - \lambda^*\|_2 \bar{\Delta}_2 
    \sum_{p=1}^k \bar{\Delta}_{1_p}
  \nonumber \\
  &\leq
    (C_t + \|\lambda - \lambda^*\|_2
    \bar{\Delta}_2) \Big(
      \frac{2}{\sigma^2} \|a\|_2 \|t(\lambda)-t^*\|_2
      + \frac{\sqrt{2}}{\sqrt{e}\sigma} \|\nu-a\|_1 
      \nonumber \\
      &\quad + \frac{2}{\sigma} \|\nu-a\|_2 \|t(\lambda)-t^*\|_2
    \Big)
    + \frac{\sqrt{2k}}{\sqrt{e}\sigma} \|a\|_2 \bar{\Delta}_2
  \nonumber \\
  &\leq
    (C_t + \|\lambda - \lambda^*\|_2 \bar{\Delta}_2) \Big(
      \frac{2\sqrt{k} C_{t^*}}{\sigma^2} \|a\|_2 \|\lambda-\lambda^*\|_2
      + \frac{\sqrt{2}}{\sqrt{e}\sigma} \|\nu-a\|_1
      \nonumber \\
      &\quad + \frac{2\sqrt{k}C_{t^*}}{\sigma} \|\nu-a\|_2 \|\lambda-\lambda^*   \|_2
    \Big)
    + \frac{\sqrt{2k}}{\sqrt{e}\sigma} \|a\|_2 
    \|\lambda - \lambda^*\|_2 \bar{\Delta}_2
    =: \bar{E}_1
    \label{eq:E1 bnd}
\end{align}
for $i=1,\ldots,k$ and $j=1,\ldots,2k$.

\subsection*{Entries in $J_2$}

By adding and subtracting $t_j^*$ then taking a Taylor
expansion like before, we obtain:
\begin{align}
  J_{2_{j,i}} &= -\phi(t_i^* -s_j + t_i(\lambda)- t_i^*)
  \nonumber \\ &=
  -\phi(t_i^* - s_j) - (t_i(\lambda)-t_i^*)\phi'(\xi_j)
  \nonumber \\ &=
  -\phi(t_i^* - s_j) - E_{2_{i,j}},
\end{align}
for some $\xi_j \in [t_i^*-s_j - |t_i(\lambda)-t_i^*|, t_i^*-s_j + |t_i(\lambda)-t_i^*|]$
and $E_{2_{j,i}}$ is the perturbation term. Then:
\begin{equation}
  |E_{2_{j,i}}| \leq 
  |t_i(\lambda) - t_i^*| \cdot \frac{\sqrt{2}}{\sigma\sqrt{e}},
  \label{eq:E2 bnd}
\end{equation}
for $i=1,\ldots,k$ and $j=1,\ldots,2k$.

\subsection*{Putting everything together}

We have that
\begin{align}
  \| E \|_F &=
    \sqrt{\sum_{i=1}^{k} \sum_{j=1}^{2k} E_{1_{j,i}}^2 
    + \sum_{i=1}^k \sum_{j=1}^{2k} E_{2_{j,i}}^2}
  \nonumber \\
  &\leq
    \sqrt{
      2k^2 \bar{E}_{1}^2
      + \frac{4k}{\sigma^2e} \sum_{i=1}^k 
        |t_i(\lambda)-t_i^*|^2
    }
  \nonumber \\
  &\leq
  k \bar{E}_1 \sqrt{2}
  + \frac{2\sqrt{k}}{\sigma\sqrt{e}} 
    \|t(\lambda) - t^*\|_2
  \nonumber \\
  &\leq
  k \bar{E}_1 \sqrt{2}
  + \frac{2 C_{t^*} k}{\sigma\sqrt{e}} 
    \|\lambda - \lambda^*\|_2,
\end{align}  
where we have used the bounds 
on the entries of $E_1$ and $E_2$
from \eqref{eq:E1 bnd} and \eqref{eq:E2 bnd}
and Theorem~\ref{thm:t dep lambda}, so this result holds 
for $\lambda \in \mathcal{B}(\lambda^*,\delta_{\lambda})$ 
for $\delta_{\lambda}$ defined in the theorem.
Finally, by substituting the expression of $\bar{E}_1$ 
from \eqref{eq:E1 bnd},
we obtain:
\begin{align}
  \|E\|_F &\leq 
 \sqrt{2} k
    \Bigg[
      \bigg(
        C_{t^*} + \|\lambda - \lambda^*\|_2 \bar{\Delta}_2
      \bigg)
      \bigg(
        \frac{2C_{t^*}\sqrt{k}}{\sigma^2} \|a\|_2 \|\lambda-\lambda^*\|_2
    \nonumber \\
        &+ \frac{\sqrt{2}}{\sqrt{e}\sigma} \|\nu-a\|_1
        + \frac{2C_{t^*}\sqrt{k}}{\sigma} \|\nu-a\|_2 \|\lambda-\lambda^*   \|_2
        \bigg)
        + \frac{\sqrt{2k}}{\sqrt{e}\sigma} \|a\|_2   
        \|\lambda - \lambda^*\|_2 \bar{\Delta}_2
    \Bigg]
    \nonumber \\
    &+ \frac{2 C_{t^*} k}{\sigma\sqrt{e}} \|\lambda-\lambda^*\|_2.
    \label{eq:E F bnd inter}
\end{align}
Let $\delta_{\gamma}$ be a bound on the perturbation:
\begin{equation}
  \left\|
    \begin{bmatrix}
      \lambda - \lambda^* \\
      \nu - a
    \end{bmatrix}
  \right\|_2
  \leq \delta_{\gamma},
\end{equation}
and therefore:
\begin{equation}
  \|\lambda-\lambda^*\|_2 \leq \delta_{\gamma}
  \qquad \text{and} \qquad
  \|\nu - a\|_2 \leq \delta_{\gamma}.
\end{equation}
We also have that:
\begin{align}
  \|\nu-a\|_2 &\leq \|\nu - a\|_1 
    \leq \|\nu\|_1 + \|a\|_1 
    \leq 2\Pi,
\end{align}
where we used that $\nu_1 +\ldots+\nu_k \leq \Pi$ 
and the 
fact that $x=\sum_{p=1}^k a_p \delta_{t_p}$ is the solution
to \eqref{eq:primal noisy}, so it 
satisfies $\|x\|_{TV} = \|a\|_1 \leq \Pi$.

Similarly, we have that:
\begin{align}
  \| \lambda - \lambda^*\|_2 &\leq
  \|\lambda\|_2 + \|\lambda^*\|_2 
  \leq
  \sqrt{k} \|\lambda\|_{\infty}
  + \sqrt{k} \|\lambda^*\|_{\infty} 
  \nonumber \\
  &\leq
  2 \sqrt{k} \tau,
\end{align}
since both $\lambda$ and $\lambda^*$ satisfy
the constraint in \eqref{eq:dual noisy}.
In order to write the bound \eqref{eq:E F bnd inter} 
as $P \cdot \delta_{\gamma}$, we expand the parentheses
and use the following bounds:
\begin{align}
  &\|\lambda-\lambda^*\|_2 \|\nu-a\|_1 
    \leq 2\Pi \cdot \delta_{\gamma} 
  \\
  &\|\lambda-\lambda^*\|_2^2 \|\nu-a\|_2 
    \leq 4\sqrt{k}\tau\Pi \cdot \delta_{\gamma}
  \\
  &\|\lambda-\lambda^*\|_2 \|\nu-a\|_2
    \leq 2\Pi \cdot \delta_{\gamma}
\end{align}
to obtain:
\begin{align}
  \|E\|_F &\leq 
    \sqrt{2} k \bigg(
      \frac{2 \sqrt{k} C_{t^*}^2 \Pi}{\sigma^2}
      + \frac{\sqrt{2k} C_{t^*}}{\sqrt{e}\sigma}
      + \frac{4 \sqrt{k} C_{t^*}^2 \Pi}{\sigma}
  \nonumber \\ 
      &+ \frac{4 k C_{t^*} \bar{\Delta}_2 \tau \Pi}{\sigma^2}
      + \frac{2\sqrt{2} \bar{\Delta}_2 \Pi}{\sqrt{e} \sigma}
      + \frac{8 k C_{t^*} \bar{\Delta}_2 \tau \Pi}{\sigma}
      + \frac{\sqrt{2k} \bar{\Delta}_2 \Pi}{\sqrt{e}\sigma}
  \nonumber \\
    &+ \frac{\sqrt{2} C_{t^*}}{\sigma\sqrt{e}}
  \bigg) \cdot \delta_{\gamma},
  \label{eq:bnd P inter}
\end{align}
which we rearrange based on $\sigma$ 
to obtain $\|E\|_F \leq P(k,\sigma,\Pi,\tau,C_{t^*}) \cdot \delta_{\gamma}$,
where:
\begin{align*}
  P(k,\sigma,\Pi,\tau,&C_{t^*}) =
  \sqrt{2} k
  \Bigg[
    \frac{1}{\sigma^2} 
    \left(
      2 \sqrt{k} C_{t^*}^2 \Pi
      + 4 k C_{t^*} \bar{\Delta}_2 \tau \Pi
    \right)
  \nonumber \\
  &+ \frac{1}{\sigma} 
  \left(
    \frac{\sqrt{2k} C_{t^*}}{\sqrt{e}}
    + 4 \sqrt{k} C_{t^*}^2 \Pi
    + \frac{2\sqrt{2} \bar{\Delta}_2 \Pi}{\sqrt{e}}
    + 8 k C_{t^*} \bar{\Delta}_2 \tau \Pi
    + \frac{\sqrt{2k} \bar{\Delta}_2 \Pi}{\sqrt{e}}
    + \sqrt{\frac2e} C_{t^*}
  \right)
  \Bigg],
\end{align*}
which is the final bound in \eqref{eq:E_F bnd}.

\section{Numerical experiments}
\label{sec:lvl numerics}

In this section, we present numerical experiments which verify 
the bounds given by our main results,  
Theorem~\ref{thm:t dep lambda}, Theorem~\ref{thm:a dep t}
and Theorem~\ref{thm:pert noisy}.
To do this, we take an example of a source and sample configuration and
a Gaussian kernel for a given $\sigma$ and solve the 
exact penalty formulation \eqref{eq:exact pen}
of the dual problem \eqref{eq:dual noisy} using
the level method~\cite{NesterovOpt}, given 
in Appendix~\ref{apdx:level method}.
We introduce
inaccuracies in $\lambda$ by stopping the algorithm early
and show how these perturbations affect the source locations
and weights. Next, we add noise to the measurements to
show how $\lambda$ is affected. We are, therefore, 
able to compare the 
ratios of the perturbations obtained numerically with the 
constants in the theorems to show the validity of our 
results in practice. The specific details are discussed 
in the next subsections.

\subsubsection*{Setup}

\begin{figure}[H]
  \centering
  \includegraphics[width=0.45\textwidth]{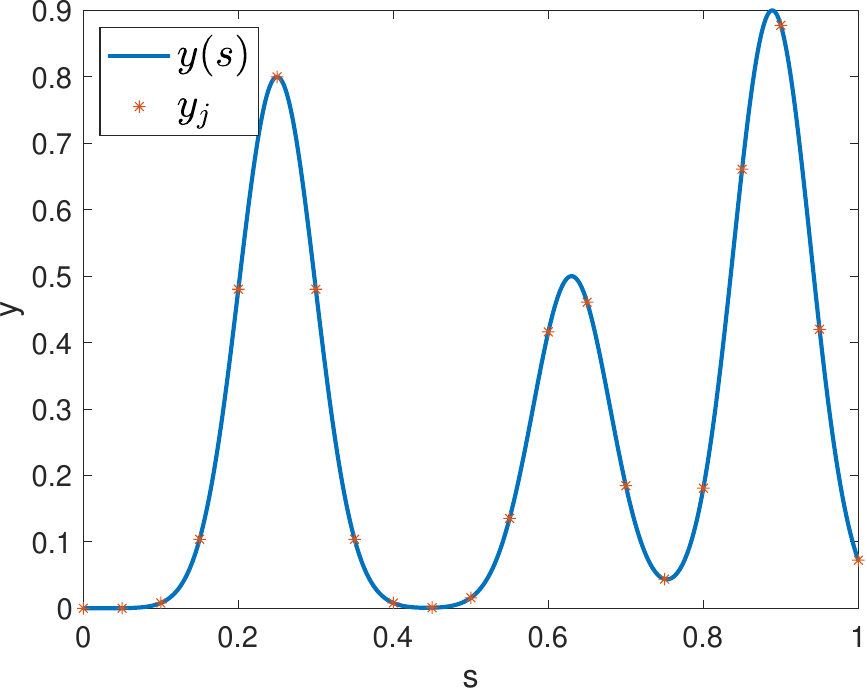}
  \caption{The source-sample configuration used
    for numerical experiments in the current section.}
  \label{fig:pert numerics cfg}
\end{figure}

We place three sources at locations $t_i^* \in T = \{0.25, 0.63, 0.889\}$ 
with weights $a_i^*~\in~\{0.8, 0.5, 0.9\}$ 
and $m=21$ equispaced samples in $[0,1]$,
with a Gaussian kernel $\phi(t) = e^{-t^2/\sigma^2}$ 
with $\sigma = 0.07$.
We show this configuration in Figure \ref{fig:pert numerics cfg}.

\subsubsection*{Effect of $\lambda^*$ perturbations on $t^*$}

We then solve the dual problem \eqref{eq:dual noisy} 
in the exact penalty formulation \eqref{eq:exact pen}
with box constraint parameter $\tau = 10^5$ and
penalty parameter $\Pi = 100$ and run
it for $P=500$ iterations. This gives 
an accuracy in the source locations of $|t_i-t^*_i| \leq 10^{-8}$ 
for $t^*_i \in T$. 

While it is possible to optimise
the parameters $\tau$, $\Pi$ and $P$
in order to obtain better accuracy
in the source locations $t_i$ and weights $a_i$,
it is not the aim of this section.
Note that Theorem~\ref{thm:t dep lambda} gives 
the result \eqref{eq:t dep lambda} in the form
\begin{equation*}
  |t_i - t^*_i| \leq C_{t^*} \|\lambda - \lambda^*\|_2, 
\end{equation*}
where $t^*_i \in T$ is an arbitrary true source location,
$\lambda^*$ is the solution to the dual 
problem \eqref{eq:dual noisy}\footnote{Note that the 
analysis of the dual problem \eqref{eq:dual}
from Section~\ref{sec:pert noise-free}
applies to the dual problem \eqref{eq:dual noisy}
considered in Section~\ref{sec:pert noisy}
as well, as the only difference difference between
\eqref{eq:dual} and \eqref{eq:dual noisy} is a box
constraint on $\lambda$.} 
and $t$ is obtained by perturbing $t^*$ as a consequence
of the perturbation $\lambda^*$ in $\lambda$.

One way of showing that a relationship of the 
type of \eqref{eq:t dep lambda} holds in practice
is to plot the ratio 
$\frac{|t^{(p)}_i - t_i^*|}{\|\lambda^{(p)}-\lambda^*\|_2}$,
for $p = p_0,\ldots,P$ and $i=1,\ldots,k$, where $P$ 
is the number of iterations the level 
method is run for, $p$ is the index of each
iteration and $t^{(p)}_i$ and $\lambda^{(p)}$ are the
values of $t_i$ and $\lambda$ obtained at iteration $p$,
where $p_0 \geq 1$ is large enough so 
that $\|\lambda^{(p)} - \lambda^*\|_2$ satisfies 
the condition in Theorem~\ref{thm:t dep lambda}.
The level method computes the value $\lambda^{(p)}$
after $p$ iterations and $\{t^{(p)}_i\}_{i=1}^k$ 
are obtained by calculating the global maxima of the 
dual certificate $q^{(p)}(s) = \sum_{j=1}^m \lambda_j^{(p)} \phi(s - s_j)$. 
Since we know the true value of $t^*_i$, we can find $t_i^{(p)}$
by running a local optimisation algorithm with $t_i^*$
as the initial condition. For a large enough value
of $p$, this will give an accurate value of $t_i^{(p)}$
and we can, therefore, calculate $|t_i^{(p)}-t_i^*|$ for
each $p = p_0,\ldots,P$ and $t_i^* \in T$.
Then we check that:
\begin{equation}
  \frac{|t_i^{(p)} - t_i^*|}{\|\lambda^{(p)}-\lambda^*\|_2}
  \leq C_{t^*},
\end{equation}
for $p = p_0,\ldots,P$ and $i=1,\ldots,k$.
One issue is that the true value of $\lambda^*$ 
is not known. The best estimate we have 
is $\lambda_{best}^* = \lambda^{(P)}$, namely the value
of $\lambda^*$ given by the level method 
after $P$ iterations. Therefore, the result 
of Theorem~\ref{thm:t dep lambda} cannot be verified 
directly in practice, but must be adapted to 
take into account this inaccuracy. 
For $i=1,\ldots,k$, we have that:
\begin{align}
  |t_i^{(p)} - t_i^*| 
  &\leq C_{t^*} \|\lambda^{(p)}-\lambda^*\|_2
  \nonumber \\
  &\leq C_{t^*} \left( \|\lambda^{(p)}-\lambda^*_{best} \|_2
    + \| \lambda^*_{best} - \lambda^*\|_2 \right),
\end{align}
and so
\begin{align}
  \frac{|t_i^{(p)} - t_i^*|}{ \|\lambda^{(p)}-\lambda^*_{best} \|_2 }
  &\leq C_{t^*} \left( 1  
    + \frac{\| \lambda^*_{best} - \lambda^*\|_2 }
    {\|\lambda^{(p)}-\lambda^*_{best} \|_2} \right).
\end{align}
For fixed $P$, which in the experiments in this section 
is $P=500$, $\| \lambda^*_{best} - \lambda^*\|_2$
above is fixed and as $p$ approaches $P$, we have
that $\| \lambda^{(p)} - \lambda^*_{best}\|_2 \to 0$, 
and therefore the right hand side above goes to infinity.
This is not a problem for our results, as it is not relevant
how the ratio 
$\frac{|t_i^{(p)} - t_i^*|}{ \|\lambda^{(p)}-\lambda^*_{best} \|_2}$
behaves for 
$\|\lambda^{(p)} - \lambda^*_{best}\|_2 \leq \|\lambda^*_{best} - \lambda^*\|_2$.

We can then find a range for $p$ where 
$\frac{\| \lambda^*_{best} - \lambda^*\|_2 }
{\|\lambda^{(p)}-\lambda^*_{best} \|_2} \leq 1$ and where we can see
that 
\begin{align}
  \frac{|t_i^{(p)} - t_i^*|}{ \|\lambda^{(p)}-\lambda^*_{best} \|_2 }
  &\leq 2 C_{t^*}.
\end{align}

\begin{figure}[!htbp]
  \centering
  \begin{minipage}[b]{0.45\linewidth}
    \centering
    \centerline{\includegraphics[width=0.93\textwidth]
      {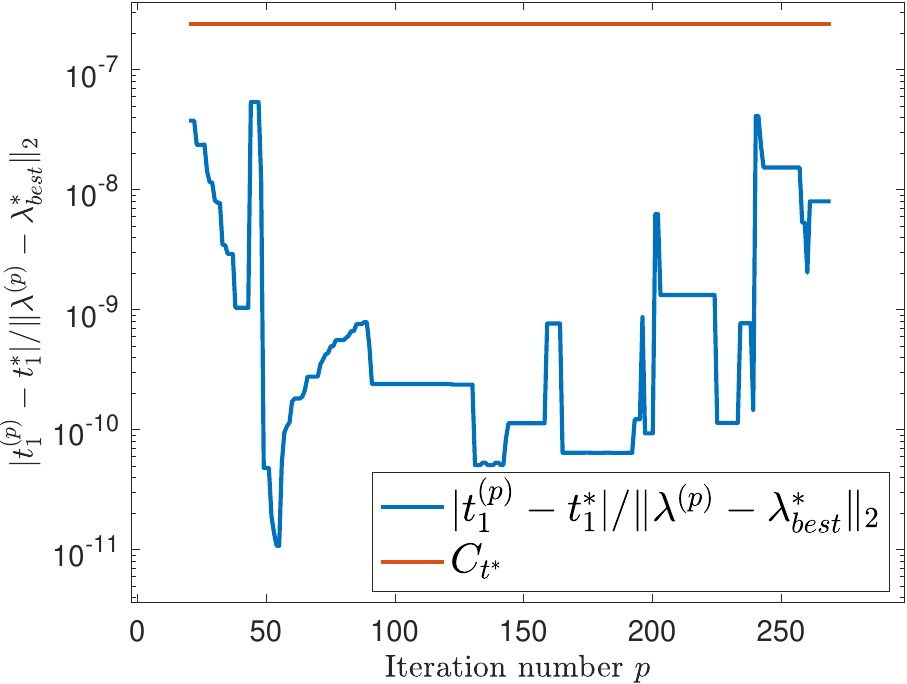}}
    \centerline{(a)}\medskip
  \end{minipage}
  \begin{minipage}[b]{0.45\linewidth}
    \centering
    \centerline{\includegraphics[width=0.93\textwidth]
      {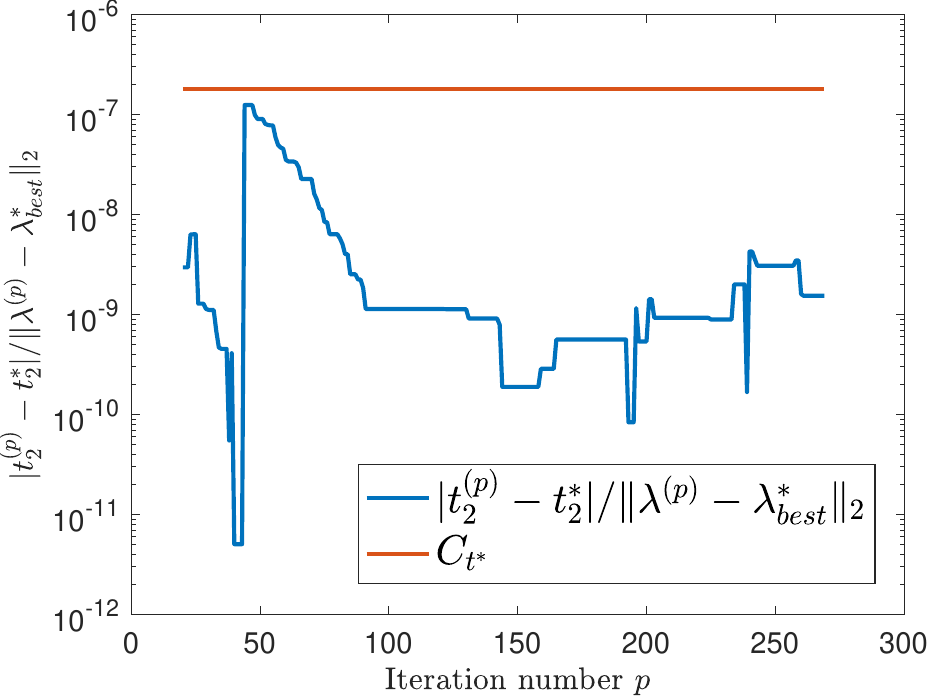}}
    \centerline{(b)}\medskip
  \end{minipage}
  \begin{minipage}[b]{0.45\linewidth}
    \centering
    \centerline{\includegraphics[width=0.93\textwidth]
      {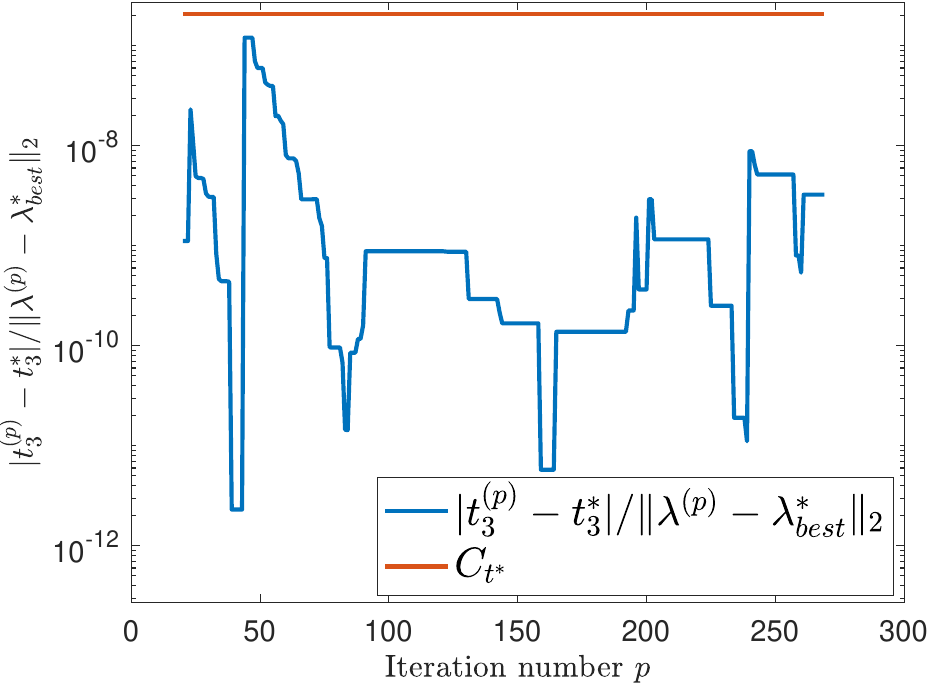}}
    \centerline{(c)}\medskip
  \end{minipage}
  \caption{The result of Theorem~\ref{thm:t dep lambda}
    for $ T = \{0.25,0.63, 0.888\}$, $\sigma=0.07$ and $m=21$.
    For each $i\in \{1,2,3\}$, we show the ratio of the 
    error in $t_i$ and the error in $\lambda$ compared
    to the constant $C_{t^*}$
    given by Theorem~\ref{thm:t dep lambda}.
  }
  \label{fig:pert t pert lam base}
\end{figure}

In Figure \ref{fig:pert t pert lam base}, we plot 
$\frac{|t_i^{(p)} - t_i^*|}{ \|\lambda^{(p)}-\lambda^*_{best} \|_2 }$
for $p = 20,\ldots,270$, where we see that the ratio is
less than $C_{t^*}$.
Specifically, we show the ratio
$\frac{\|t_i^{(p)} - t_i^*\|}{\|\lambda^{(p)}-\lambda_{best}^*\|_2}$
and the constant $C_{t^*}$ from Theorem~\ref{thm:t dep lambda}
for each $i \in \{1, 2, 3\}$.

\subsubsection*{Effect of $t^*$ perturbations on $a^*$}

\begin{figure}[!htbp]
  \centering
  \includegraphics[width=0.5\textwidth]
      {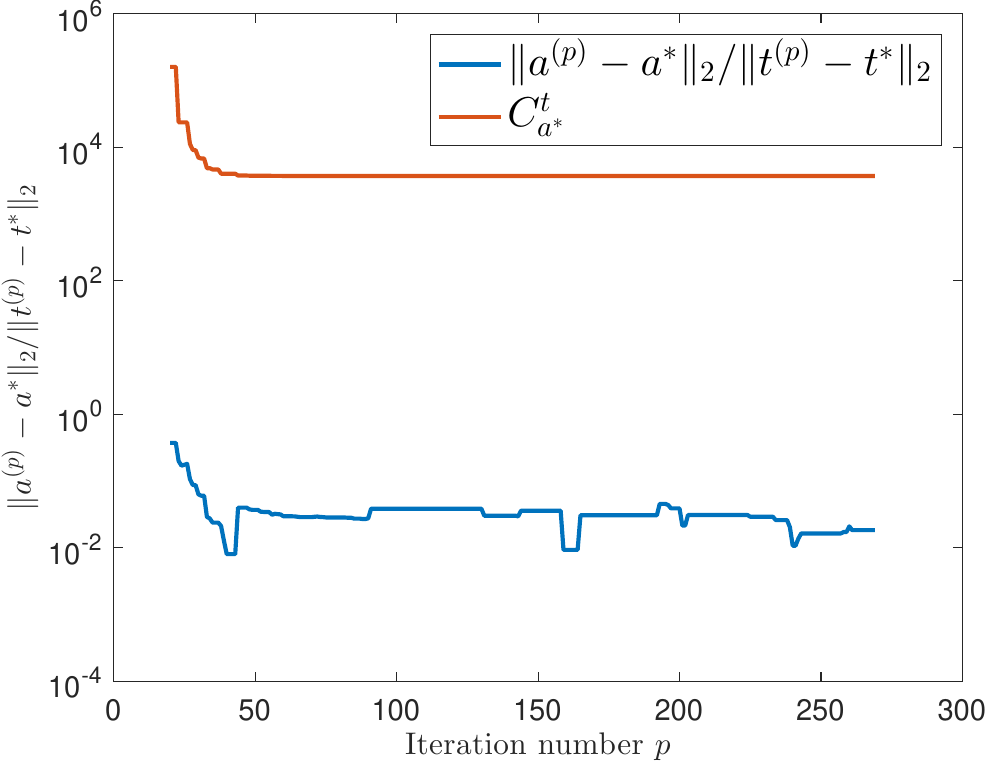}
  \caption{
    Plot of the ratio between $\|a^{(p)}-a^*\|_2$ and $\|t^{(p)}-t^*\|_2$ 
    for $p = 20,\ldots,270$, and the bound $C_{a^*}^t$  
    from \eqref{eq:ratio err a t}
    in the setup described at the beginning of this section.
  }
  \label{fig:pert a pert t base}
\end{figure}

In the case of Theorem~\ref{thm:a dep t}, it is more straightforward
to check the ratio of the errors, since we know the true
values of the source locations and weights, which we denote 
by $t^* = [t_1^*,\ldots,t_k^*]^T$ and $a^* = [a_1^*,\ldots,a_k^*]^T$ respectively.
The error bound \eqref{eq:a dep t} given by the theorem is of the form:
\begin{equation*}
  \| a - a^*\|_2 \leq C_{a^*} 
  e^{\frac{4\|t-t^*\|_{\infty}}{\sigma^2}}
  \|t - t^*\|_2 + \mathcal{O}(\|t - t^*\|_2^2),
\end{equation*}
where $t$ is the perturbed vector $t^*$ and $a$ is the perturbed
vector $a^*$ as a consequence of perturbing $t^*$.
For the values $t^{(p)}_i, i \in \{1,2,3\}$, obtained after $p$ 
iterations of the level method, we now solve the least squares 
problem $argmin_{\hat{a}} \|\Phi^{(p)}\hat{a}-y\|_2$ 
with the entries in the data matrix $\Phi^{(p)}$ 
given by $\Phi^{(p)}_{j,i} = \phi(t^{(p)}_i-s_j)$ to find
the corresponding perturbed weights $a^{(p)}_i$ for $i \in \{1,2,3\}$.
Then, according to Theorem~\ref{thm:a dep t}, we have that:
\begin{equation}
  \frac{\|a^{(p)} - a^*\|_2}{\|t^{(p)} - t^*\|_2} 
  \leq C_{a^*}^t + \mathcal{O}(\|t^{(p)} - t^*\|_2),
  \label{eq:ratio err a t}
\end{equation}
where we write
\begin{equation*}
  C_{a^*}^t = C_{a^*} e^{\frac{4\|t^{(p)}-t^*\|_{\infty}}{\sigma^2}}.
\end{equation*}
In Figure \ref{fig:pert a pert t base}, we show
the ratio $\frac{\|a^{(p)} - a^*\|_2}{\|t^{(p)}-t^*\|_2}$
and $C_{a^*}^t$ in
the same setting as in Figure \ref{fig:pert t pert lam base},
for iterations $p=20,\ldots,270$.

\subsubsection*{Effect of the noise $w$ on $\lambda^*$ and $t^*$}

As in the case of Theorem~\ref{thm:t dep lambda}, where we rely on a best
approximation $\lambda^*_{best}$ of $\lambda^*$ for the 
numerical experiments, a similar approach is required
to check the validity of the results of Theorem~\ref{thm:pert noisy} 
in practice. Theorem~\ref{thm:pert noisy} gives the 
bound \eqref{eq:pert noisy} in the form:
\begin{equation*}
  \| \lambda^*_{w} - \lambda^*\|_2 
  \leq C_{\lambda^*} \cdot \|w\|_2,
\end{equation*}
where $\lambda^*$ is the true solution of the dual
problem \eqref{eq:dual noisy} and $\lambda^*_{w}$ 
is the solution to the same problem
with $y$ perturbed by the noise $w$.

\begin{figure}[!htbp]
  \centering
  \begin{minipage}[b]{0.45\linewidth}
    \centering
    \centerline{\includegraphics[width=1\textwidth]
      {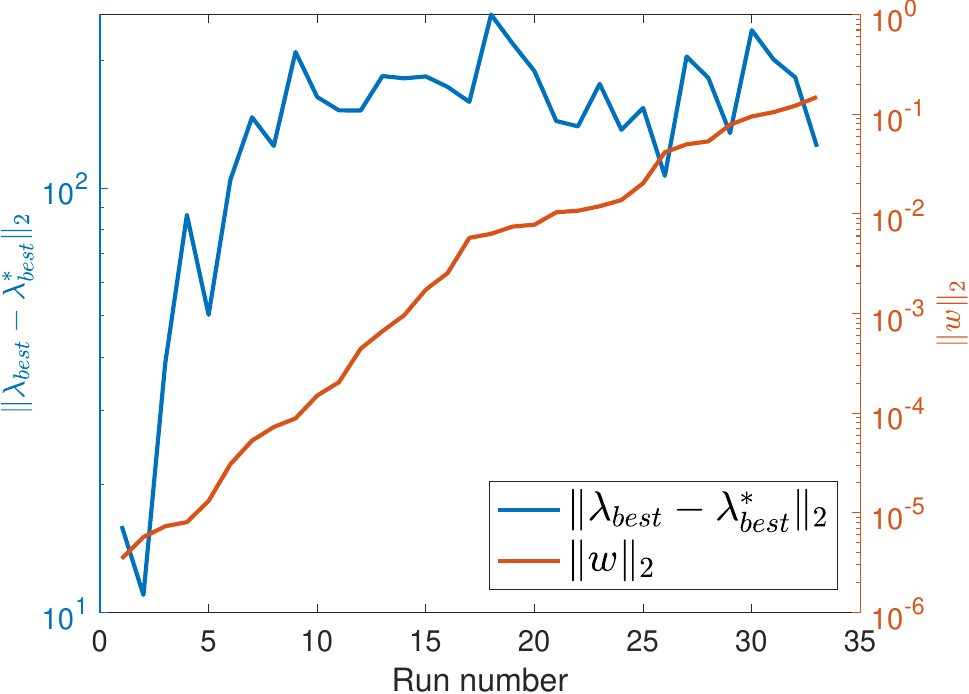}}
    \centerline{(a)}\medskip
  \end{minipage}
  \hspace{0.5cm}
  \begin{minipage}[b]{0.45\linewidth}
    \centering
    \centerline{\includegraphics[width=1\textwidth]
      {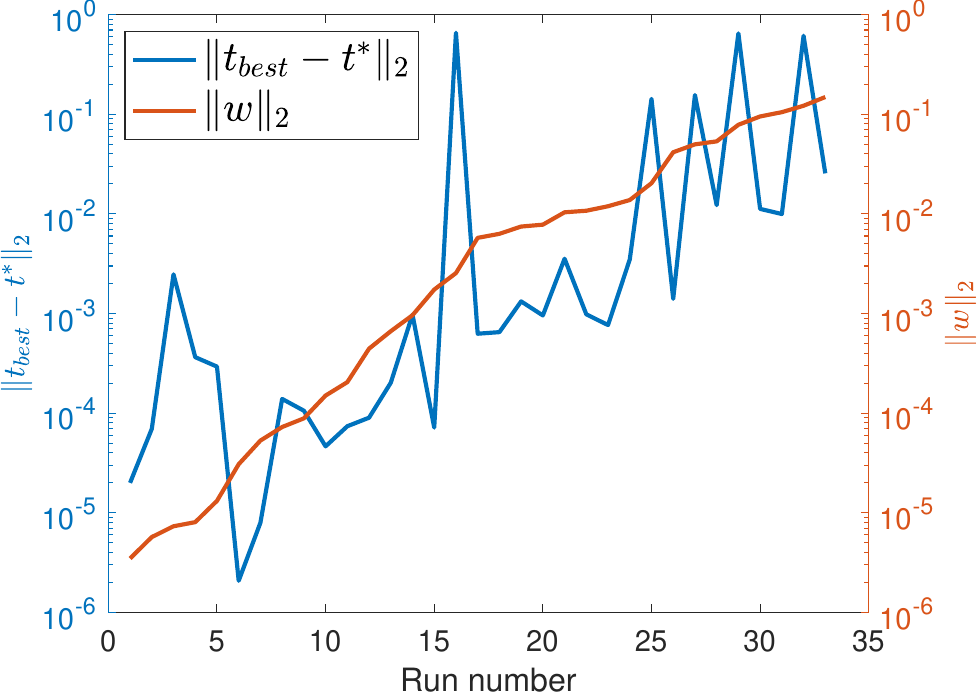}}
    \centerline{(b)}\medskip
  \end{minipage}
  \begin{minipage}[b]{0.45\linewidth}
    \centering
    \centerline{\includegraphics[width=1\textwidth]
      {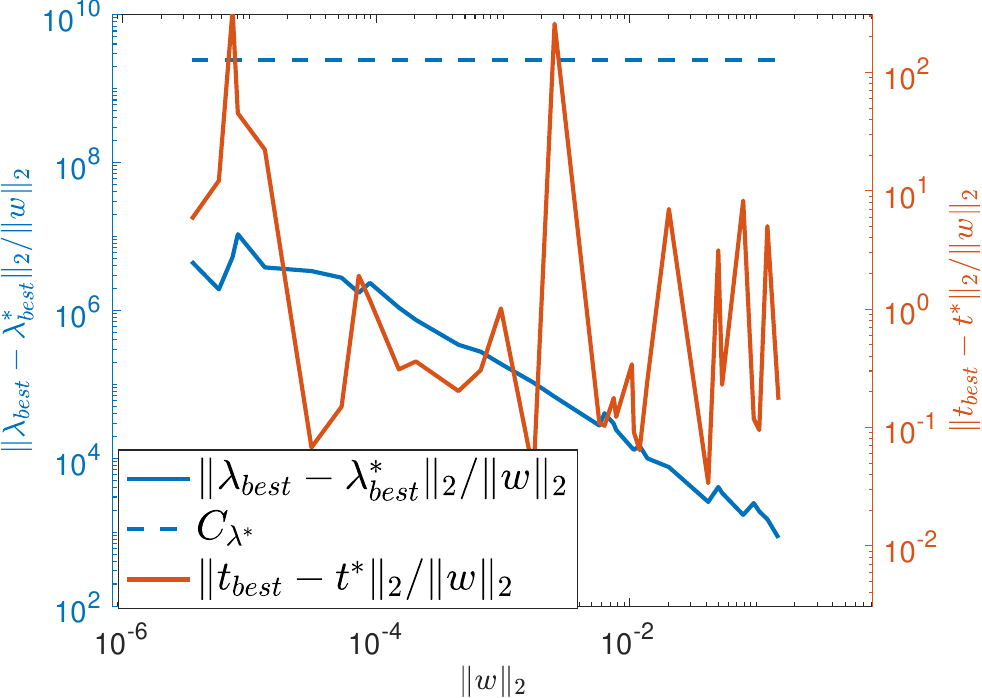}}
    \centerline{(c)}\medskip
  \end{minipage}
  %
  \caption{Plots of $\|\lambda_{best}-\lambda^*_{best}\|_2$
    (panel(a)), $\|t_{best}-t\|_2$ (panel (b)) 
    and their ratio to the noise $\|w\|_2$ (panel(c))
    for $\|w\|_2$ in a range as given in 
    \eqref{eq:calc noise} and \eqref{eq:noise range}, 
    in the setting described at the beginning of this section.}
  \label{fig:t lam pert w base}
\end{figure}

As it is not possible to know exactly the values 
of $\lambda^*$ and $\lambda^*_{w}$, 
let $\lambda^*_{best} = \lambda^{(P)}$ be the value
of $\lambda$ given by the level method after $P$ iterations
when $y$ is exact and $\lambda_{best}$ be the value
of $\lambda$ returned by the level method after $P$
iterations when $y$ is corrupted by the additive noise $w$.
Then we can reformulate the bound \eqref{eq:pert noisy}
in terms of $\lambda^*_{best}$ and $\lambda_{best}$:
\begin{align}
  \|\lambda_{best}-\lambda^*_{best}\|_2 
  &= 
  \|\lambda_{best}-\lambda^*_{best}
    + \lambda^* - \lambda^* 
    + \lambda^*_{w} - \lambda^*_{w} 
  \|_2 
  \nonumber \\
  &\leq
  \|\lambda_{best} - \lambda^*_{w}\|_2
  + \|\lambda^* - \lambda^*_{w}\|_2
  + \|\lambda^* - \lambda^*_{best}\|_2
  \nonumber \\
  &\leq
  \|\lambda_{best} - \lambda^*_{w}\|_2
  + C_{\lambda^*} \|w\|_2
  + \|\lambda^* - \lambda^*_{best}\|_2,
\end{align}
so 
\begin{equation}
  \frac{
    \|\lambda_{best}-\lambda^*_{best}\|_2 
  }{\|w\|_2}
  \leq
  C_{\lambda^*} +
  \frac{
    \|\lambda_{best} - \lambda^*_{w}\|_2
    + \|\lambda^* - \lambda^*_{best}\|_2
  }{\|w\|_2}.
  \label{eq:Clam plus inf}
\end{equation}
As before, we plot 
$  
  \frac{
    \|\lambda_{best}-\lambda^*_{best}\|_2 
  }{\|w\|_2}
$,
where $\lambda^*_{best}$ is the solution we obtain by
solving the dual problem \eqref{eq:dual noisy}
in its exact penalty formulation using the level 
method with $P=100$ iterations and $\lambda_{best}$
is the `noisy' solution, which is obtained
by solving the problem with $P=100$ iterations
when $y$ is corrupted by additive noise $w$.
We repeat this for different magnitudes of the
noise $w$, which we increase gradually as 
follows. For each
component $y_j$ of $y$, we add a sample $X_j$
from the standard uniform distribution $U(0,1)$,
multiplied by a coefficient $w_c$:
\begin{equation}
  y_{{noisy}_j} = y_j + w_c \cdot X_j.
  \label{eq:calc noise}
\end{equation}
We repeat this for different values 
of the coefficient $w_c$ from the set:
\begin{align}
  w_c \in \{
    &0.000002, 0.000004,\ldots,0.00001,
    \nonumber \\
    &0.00002, 0.00004,\ldots,0.0001,
    \nonumber \\
    &0.0002, 0.0004,\ldots,0.001,
    \nonumber \\
    &0.002, 0.003,\ldots,0.01,
    \nonumber \\
    &0.02, 0.03,\ldots,0.1
  \}.
  \label{eq:noise range}
\end{align}
Therefore, in Figure \ref{fig:t lam pert w base} 
we show the basic setup described at the beginning of 
this section. 
Panel (a) shows $\|\lambda_{best} - \lambda_{best}^*\|_2$
against the norm of the noise $\|w\|_2$, and in order
to check that the algorithm actually converges to a
useful $\lambda_{best}^*$, we also plot plot
$\|t_{best} - t^*\|_2$ against $\|w\|_2$
in panel (b), since we know the true value $t^*$.
Then, in panel (c) we plot the ratio
$\frac{\|\lambda_{best} - \lambda_{best}^*\|_2}{\|w\|_2}$
and $C_{\lambda^*}$ as given by Theorem~\ref{thm:pert noisy},
where we see that the ratio is smaller than the constant, 
as the theorem states. In the same plot, we also show 
the ratio $\frac{\|t_{best}-t^*\|_2}{\|w\|_2}$
and we see that it does not grow as the magnitude
of the noise increases. 
In these experiments
we only take into account $2k$ entries 
of $\lambda$ and $w$, corresponding to the $2k$ samples
that are the closest to the $k$ sources,
as described in Section~\ref{sec:pert noisy}, 
for which Theorem~\ref{thm:pert noisy} holds.

\section{Conclusion}
\label{sec:conclusion}

We have developed a stability analysis of the dual approach to non-negative super-resolution with a Gaussian kernel. 
Using a quantitative version of the implicit function theorem, we tracked how error propagates along
the recovery pipeline: 
from measurement noise, to the solution of the convex dual, to the recovered spike locations, and finally to their weights. 
In the noise-free setting, we bounded the accuracy of the primal solution in terms of
the accuracy of the dual solution and the number of measurements, with an
$\ell_\infty$ bound on the spike locations (Theorem~\ref{thm:t dep lambda}) and
an $\ell_2$ bound on the weights (Theorem~\ref{thm:a dep t}). 
In the noisy setting, we bounded the perturbation of the dual variable as a function of the
noise level (Theorem~\ref{thm:pert noisy}), closing the chain. The constants are
explicit throughout, and the curvature $q''(t^*)$ of the dual certificate at the
sources emerges as the quantity that governs how well the locations can be
resolved.

The noisy analysis is made tractable by the $\ell_1$ data fidelity formulation:
its convex dual is box-constrained, and the complementarity at optimality that this induces is exactly what yields the square, invertible Jacobian on which our perturbation bound rests. 
Several questions remain open. 
The dependence of our bounds on the number of measurements stems from a deliberately loose control of sums of shifted Gaussians and is likely not tight, and we do not establish a
dependence on the separation between the sources. 
Finally, since only kernel-specific bounds enter the analysis, the same strategy should extend to other kernels.

\section*{Acknowledgements}

This work was done while BT was affiliated to the
Mathematical Institute, University of Oxford, UK.
This publication is based on work supported by the 
EPSRC Centre For Doctoral Training in Industrially Focused 
Mathematical Modelling (EP/L015803/1) in collaboration with 
the National Physical Laboratory and by the Alan Turing Institute 
under the EPSRC grant EP/N510129/1 and the Turing Seed Funding grant SF019. 

\bibliographystyle{unsrt}
\bibliography{main.bbl}

\appendix

\section{Duality in the noisy case}
\label{apdx:dual noisy} 

In this section, we show the duality of problems~\eqref{eq:primal noisy} and~\eqref{eq:dual noisy raw}, namely
\begin{equation*}
  \min_{x \geq 0} 
  \left\|
    y - \int \Phi(t) x(\dif t)
  \right\|_1
  \quad\text{subject to}\quad
  \| x \|_{TV} \leq \Pi,
\end{equation*}
and
\begin{equation*}
  \max_{\substack{\beta > 0\\\lambda \in \mathbb{R}^m}}
  \beta \left(
    \lambda^T y
    - \Pi
  \right)
  \quad\text{subject to}\quad
  \lambda^T \Phi(t) \leq 1,
  \quad\forall t \in [0,1]
  \quad\text{and}\quad 
  \|\lambda\|_{\infty} \leq 1/\beta.
\end{equation*}
We start by introducing a new variable $z=\int \Phi(t)x(\dif t)$ in~\eqref{eq:primal noisy}:
\begin{align}
  \min_{\substack{x \geq 0\\z \in \mathbb{R}^m}}
  \left\|
    z-y
  \right\|_1
  \quad \text{subject to} \quad
  &z = \int \Phi(t) x(\dif t),
  \nonumber \\
  &\| x \|_{TV} \leq \Pi,
\end{align}
and writing the Lagrangian:
\begin{align}
  L(x,z,\beta,\lambda) = 
  \|z-y\|_1 + \lambda^T 
  \left(
    z - \int \Phi(t) x(\dif t)
  \right)
  + \beta \left( \|x\|_{TV} - \Pi \right).
\end{align}
Then, the Lagrangian dual problem is:
\begin{align}
  \max_{\substack{\beta \geq 0\\\lambda \in \mathbb{R}^m}}
  &\min_{\substack{x \geq 0\\z \in \mathbb{R}^m}} 
  L(x,z,\beta,\lambda) = 
  \nonumber \\
  &=
  \max_{\substack{\beta \geq 0\\\lambda \in \mathbb{R}^m}}
  \min_{\substack{x \geq 0\\z \in \mathbb{R}^m}} 
  \left[ 
    \|z - y\|_1 + \lambda^T z
    + \int \left(
      \beta - \lambda^T \Phi(t) 
    \right)x(\dif t)
  \right]
  - \beta \Pi
  \nonumber \\
  &=
  \max_{\substack{\beta \geq 0\\\lambda \in \mathbb{R}^m}}
  \min_{\substack{x \geq 0\\w \in \mathbb{R}^m}} 
  \left[ 
    \|w\|_1 + \lambda^T w
    + \int \left(
      \beta - \lambda^T \Phi(t) 
    \right)x(\dif t)
  \right]
  + \lambda^T y 
  - \beta \Pi,
\end{align}
where in the last equality we make the 
substitution $w = z-y$.

The integral on the right hand side is equal 
to $-\infty$ if there exists $t_0 \in [0,1]$ such
that $\lambda^T \Phi(t_0) > \beta$, as we can 
set $x = \infty \cdot \delta_{t_0}$. Therefore, 
we impose the condition that
$\lambda^T \Phi(t) \leq \beta$ for all $t \in [0,1]$,
in which case the integral is equal to zero 
by taking $x$ to be zero wherever the integrand is non-zero,
and the dual becomes:
\begin{align}
  \max_{\substack{\beta \geq 0\\\lambda \in \mathbb{R}^m}}
  \min_{w \in \mathbb{R}^m}
  \left(
    \|w\|_1 + \lambda^T w
  \right)
  + \lambda^T y
  - \beta \Pi
  \quad\text{subject to}\quad
  \lambda^T \Phi(t) \leq \beta,
  \quad\forall t \in [0,1].
\end{align}
which can be rewritten as:
\begin{align}
  \max_{\substack{\beta \geq 0\\\lambda \in \mathbb{R}^m}}
  -\max_{w \in \mathbb{R}^m}
  \left\{
     -\lambda^T w -\|w\|_1
  \right\}
  + \lambda^T y
  - \beta \Pi
  \quad\text{subject to}\quad
  \lambda^T \Phi(t) \leq \beta,
  \quad\forall t \in [0,1].
\end{align}
and note that for $f(w) = \|w\|_1$:
\begin{equation}
  f^*(\lambda) = \max_{w} \left\{
    \lambda^T (-w) - \|-w\|_1 
  \right\}
  = \begin{cases}
    0, &\text{if} \quad \|\lambda\|_{\infty} \leq 1,
    \\
    \infty, &\text{otherwise},
  \end{cases}
\end{equation}
is its conjugate \cite{Hiriart-Urruty1996}.
Therefore, we impose the condition 
that $\|\lambda\|_{\infty} \leq 1$ and the dual becomes:
\begin{align}
  \max_{\substack{\beta \geq 0\\\lambda \in \mathbb{R}^m}}
  \lambda^T y
  - \beta \Pi
  \quad\text{subject to}\quad
  \lambda^T \Phi(t) \leq \beta,
  \quad\forall t \in [0,1]
  \quad\text{and}\quad 
  \|\lambda\|_{\infty} \leq 1.
\end{align}
We then make the substitution $\lambda' = \lambda/\beta$ 
(for $\beta > 0$)
to obtain:
\begin{align}
  \max_{\substack{\beta > 0\\\lambda' \in \mathbb{R}^m}}
  \beta \left(
    \lambda'^T y
    - \Pi
  \right)
  \quad\text{subject to}\quad
  \lambda'^T \Phi(t) \leq 1,
  \quad\forall t \in [0,1]
  \quad\text{and}\quad 
  \|\lambda'\|_{\infty} \leq 1/\beta,
\end{align}
which is the problem \eqref{eq:dual noisy raw}.
Fixing $\beta$ and keeping the argmax over $\lambda$ gives~\eqref{eq:dual noisy}, see Section~\ref{sec:main goals}.

\section{The level bundle method}
  \label{apdx:level method}

In this section, we describe the level bundle 
method~\cite{NesterovOpt} applied to \eqref{eq:exact pen} 
for which experiments were presented in 
Section~\ref{sec:main goals} and Section~\ref{sec:lvl numerics}. 
The algorithm progressively builds up a polyhedral model 
of the objective function from a `bundle' of subgradients 
at each iteration. The algorithm proceeds by projecting iterates onto a level set of the model, an approach which is known to improve robustness in comparison with the standard cutting planes subgradient method (Kelley's method). A statement of the algorithm is given in Algorithm~\ref{alg:level method alg}.

\begin{algorithm}[h]
\vspace{5pt}
\textbf{Input:} Kernel function $\Phi:I\rightarrow\mathbb{R}^{m}$, measurements $y\in\mathbb{R}^{m}$, sample locations $\{s_j\}_{j\in\{1,\ldots,m\}}\in I$, penalty parameter $\Pi>0$, level set parameter $\alpha\in(0,1)$ and number of iterations $L$.\\
 \\
\textbf{Initialize:} $l=1$. \\
\\
\textbf{While} $l\le L$ , \textbf{do}
\begin{enumerate}
\item Compute a subgradient as
$$t^l\in\arg\sup_{s\in I}\displaystyle(\lambda^{l-1})^T\Phi(s),$$
$$g^l=\left\{\begin{array}{ll}-y+\Pi
    \left[(\lambda^{l-1})^T\Phi(t^l)-1\right],&(\lambda^{l-1})^T\Phi(t^l)\geq 1\\
-y,&(\lambda^{l-1})^T\Phi(t^l)<1\end{array}\right.$$
\item Build the polyhedral model
$$\widehat{\Psi}^l_{\Pi}(\lambda)=\max_{r=1,\ldots,l}\Psi_{\Pi}(\lambda^{r-1})+(g^r)^T(\lambda-\lambda^{r-1}).$$
\item Compute $\nu^l=\displaystyle\inf_{\lambda}\widehat{\Psi}_{\Pi}^l(\lambda)$ and $\mu^l=\displaystyle\min_{r=1,\ldots,l}\Psi_{\Pi}(\lambda^{r})$.
\item Project onto the level set as $\lambda^l=\mathcal{P}_{\mathcal{L}_{\alpha}^l}(\lambda^{l-1})$ where $\mathcal{L}_{\alpha}^l=\{\lambda:\widehat{\Psi}_{\Pi}(\lambda)\le\alpha\mu^l+(1-\alpha)\nu^l\}$.
\item $l=l+1$.
\end{enumerate}
\vspace{5pt}
\textbf{Output:} $\lambda^L\in\mathbb{R}^m$.
\caption{Level bundle method for solving Program \eqref{eq:exact pen}.}
\label{alg:level method alg}
\end{algorithm}

In the experiments shown in Section~\ref{sec:main goals}, $\Pi$ was chosen to be $2\|\bf{a^{\ast}}\|_1$ and the level set parameter $\alpha$ was taken to be $1/4$.

\end{document}